\documentclass[12pt]{article}
\usepackage{verbatim,amsmath,amssymb,mathrsfs,amsthm,amsfonts,bm}
\usepackage{varwidth}
\usepackage{extarrows}
\usepackage{enumitem}
\usepackage{fancyhdr}
\usepackage{booktabs}
\usepackage[dvipsnames]{xcolor}
\usepackage{graphicx,wrapfig}
\usepackage{subfigure}
\usepackage{caption}
\usepackage{amsmath}
\usepackage{makecell}
\usepackage{multirow}
\usepackage{array}
\usepackage{setspace}
\usepackage{float} 
\usepackage{cite}
\usepackage{color,mathtools}
\usepackage{epstopdf}
\usepackage{lipsum}
\usepackage{appendix}
\usepackage{algorithm}
\usepackage{algorithmic}
\usepackage{blindtext}
\usepackage{ulem}
\allowdisplaybreaks

\numberwithin{equation}{section}
\numberwithin{figure}{section}
\numberwithin{table}{section}
\allowdisplaybreaks

\newtheorem{theorem}{Theorem}[section]
\newtheorem{lemma}[theorem]{Lemma}

\newtheorem{remark}{Remark}[section]
\newtheorem{example}{Example}[section]


\usepackage{hyperref}
\hypersetup{
    colorlinks=true, linkcolor=blue,
    urlcolor=cyan, citecolor=blue,
}

\begin{document}

\title{Structure-preserving active flux method for the compressible Euler equations with gravity}
\date{}

\author{}
\author{Jiahui Zhang\thanks{Institute for Mathematics, Johannes Gutenberg University of Mainz, Mainz, 55128, Germany.  E-mail: jizhang@uni-mainz.de.}
\and Christian Klingenberg\thanks{Corresponding author. Institute of Mathematics, University of W\"urzburg, W\"urzburg, 97074, Germany. E-mail: christian.klingenberg@uni-wuerzburg.de} 
}

\maketitle

\begin{abstract}
  In this paper, we propose a structure-preserving active flux method for the compressible Euler equations under gravitational fields, ensuring exactness for the generalized hydrostatic equilibrium state and the preservation of positivity for both density and pressure. This approach entails the update of the cell average and the point values simultaneously. To obtain the equilibrium preservation of the scheme for the cell average, the strategy pursued is to reformulate the gravitational term using the steady state solution, followed by the application of a decomposition algorithm. For the point values, we introduce minor adjustments to the Lax-Friedrichs flux vector splitting, coupled with the same finite difference operator performed on the flux and source term to restore the well-balancedness. Furthermore, a posteriori Multi-dimensional Optimal Order Detection (MOOD) paradigm is adopted for achieving positivity preservation and oscillation control near strong discontinuities. A series of numerical experiments are performed to verify the high-order accuracy of our proposed method, along with its exact equilibrium and positivity preservation, and its ability to capture small perturbations at or close to steady flow without spurious oscillations. 

  \smallskip
  \textbf{Keywords}: Euler equations with gravity; active flux method; well-balanced; positivity-preserving; source term treatment; hydrostatic equilibrium.
\end{abstract}

\section{Introduction}\label{se:in}
This paper concentrates on discussing the compressible Euler equations which incorporate the influence of the gravitational force. Many interesting phenomena in astrophysical and atmospheric are modelled by this system, such as core-collapse supernova simulation, stellar evolution, numerical weather forecasting, etc. In the $d$-dimensional case, the governing equations for the conservation of mass, momentum and energy, together with the gravitational field, are expressed as follows, 
\begin{equation}\label{euler2D}
    \left\{\begin{array}{l} 
    \rho_t+\nabla\cdot(\rho \boldsymbol{u})=0,\\
    (\rho  \boldsymbol{u})_t+ \nabla \cdot (\rho  \boldsymbol{u}\otimes \boldsymbol{u} + p \mathbf{I}_d )  =-\rho \nabla\phi,\\
    E_t+\nabla \cdot((E+p) \boldsymbol{u}) = -\rho  \boldsymbol{u}\cdot \nabla\phi,
    \end{array}\right.
\end{equation}
where $\rho, \boldsymbol{u},p$ denote the mass density, the velocity of the gas, and the pressure, respectively. $\rho \boldsymbol{u}$ corresponds to the momentum vector. The total fluid energy $E$ is defined as the sum of the internal and kinetic energy, i.e., $E=\rho e + \frac{1}{2}\rho \|\boldsymbol{u}\|^2$, where $e$ denotes the internal energy. $\phi=\phi(\boldsymbol{x})$ represents the time-independent gravitational potential. $\mathbf{I}_d$ is the identity matrix and the operators $\nabla$, $\nabla\cdot$, $\otimes$ signify the gradient, divergence and tensor product. To fully close this system of equations, it is essential to incorporate an equation of state (EoS) that establishes a relationship between the specific internal energy, density, and pressure. It is given by
\begin{equation}\label{eos}
  p = p(\rho,e) = (\gamma-1)\rho e = (\gamma-1)\left(E - \dfrac{1}{2}\rho \|\boldsymbol{u}\|^2\right),
\end{equation}
with $\gamma$ being the ratio of specific heats.

As a typical example of hyperbolic balance laws, one significant aspect of such models is their tendency to exhibit non-trivial steady state solutions, where the non-zero flux gradient precisely balances the source term. In particular, they admit hydrostatic stationary solutions with $\boldsymbol{u}= 0$, corresponding to a mechanical equilibrium where pressure balances external forces such as gravity, 
\begin{equation}\label{hydro:euler}
  \rho = \rho(\boldsymbol{x}), \quad  \boldsymbol{u}= 0, \quad \nabla p = -\rho\nabla \phi.
\end{equation}
Since this system of equations is underdetermined, additional assumptions are required to fully characterize the equilibrium. One special class of equilibrium is the isothermal hydrostatic state, where the temperature is a constant, i.e., $T(x)=T_0$. For an ideal gas, it is given by
\begin{equation}\label{iso_steady} 
  \rho = \rho_0\exp\left(-\dfrac{\phi}{RT_0}\right), \quad  \boldsymbol{u}= 0, \quad p = p_0\exp\left(-\dfrac{\phi}{RT_0}\right),
\end{equation}
where $R$ is the gas constant, $\rho_0,p_0$ are constants satisfying $p_0=\rho_0RT_0$. Another class of equilibrium state is polytropic hydrostatic state given by 
\begin{equation}\label{pol_steady}
  \rho = \left( \dfrac{\gamma-1}{K\gamma}(c-\phi) \right)^\frac{1}{\gamma-1}, \quad \boldsymbol{u}=0, \quad p = \dfrac{1}{K^{\frac{1}{\gamma-1}}}\left(\dfrac{\gamma-1}{\gamma}(c-\phi)\right)^\frac{\gamma}{\gamma-1},
\end{equation}
with $c$ being a constant, $K$ is a constant defined by the polytropic relation $p=K\rho^\gamma$.

The existence of source terms presents a well-known challenge for the standard numerical schemes, which often fail to preserve the steady state even when highly refined meshes are utilized. This promotes the development of methods that can exactly maintain steady state solutions at the discrete level while accurately resolving small perturbations on relatively coarse grids. Such methods are referred to as well-balanced schemes, first introduced by Bermudez and Vazquez \cite{bermudez1994upwind}. 
In recent years, considerable attention has been devoted to the design of well-balanced schemes for the Euler equations with gravitation. The first work was constructed by LeVeque and Bale \cite{LeVeque1998wave}, where a quasi-steady wave-propagation algorithm was extended to an ideal gas under a static gravitational field. 
Later, many efforts have been made in this direction, including the finite difference methods \cite{ghosh2016well,xing2013highfde}, finite volume methods \cite{li2016highfve,klingenberg2019arbitrary,kappeli2014well}, and discontinuous Galerkin schemes \cite{wu2021uniformly,li2018welldge,zhang2021high}. 
Another class of methods \cite{berberich2020high,varma2019second} were designed to maintain consistency between the discrete discretization and the steady-state partial differential equations of interest. It is free of additional assumptions and a priori knowledge of the hydrostatic solution. 
In addition, several studies \cite{grosheintz2020well,2023arXivequilibrium,zhang2023moving,xuhigh} have explored a generalized more complicated moving equilibrium with the non-zero velocity.

Another challenge in developing structure-preserving methods is the property of positivity preservation. Preserving the positivity of the density and pressure is essential for the robustness of numerical algorithms. Otherwise, the computation might break down, especially for low density problems such as blast waves and low pressure problems such as high Mach number astrophysical jets. In recent years, positivity-preserving techniques have gained significant popularity for high-order schemes. The most commonly used is a simple linear scaling limiter proposed by Zhang and Shu \cite{zhang2010maximum} for the reconstructed polynomials in finite volume or discontinuous Galerkin (DG) schemes. It can be generalized to the Euler equations without the source term \cite{zhang2010positivity} and with source term \cite{zhang2011positivity}. An alternative strategy is Multi-dimensional Optimal Order Detection (MOOD) paradigm, which was first introduced by Clain et al. \cite{clain2011high} and further developed by Diot et al. \cite{diot2012improved}. Instead of relying on a priori scaling of the solution, MOOD adopts a posteriori detection and recomputation strategy. When problematic solutions are detected after each time update, the scheme will reduce the local polynomial degree before recomputing the solution. 

Recently, the active flux method is a relatively new variant of a finite volume method which is attracting increasing interest. Both cell averages and point values at cell interfaces are used as degrees of freedom. This approach is characterized by a compact, globally continuous reconstruction, the avoidance of Riemann solvers at cell interfaces, while maintaining high-order accuracy and conservation properties. 
In the early developments of the active flux method, exact or approximate evolution operators are employed to evolve the point values in conjunction with Simpson's rule to compute the flux quadrature in time. Such techniques have been applied to a variety of hyperbolic conservation laws, for instance, the linear equations \cite{eymann2013multidimensional,van1977towards}, Burgers' equation \cite{barsukow2021non,chudzik2021cartesian}, one-dimensional shallow water equations \cite{barsukow2024well}, and the compressible Euler equations \cite{eymann2011active,barsukow2021non,chudzik2024active}, etc. These resulting schemes are fully discrete.
However, some limitations exist since it is difficult to obtain evolution operators for nonlinear systems of equations, especially in multi-dimensional case. This has led \cite{abgrall2023combination,abgrall2025semi,duan2025active,abgrall2023extensions} to consider semi-discrete formulations, where the cell average and point values are used to discretize the spatial derivatives and subsequently advanced in time using suitable time integration methods. As outlined in \cite{abgrall2023hybrid}, this approach applies readily to nonlinear problems free of evolution operators, at the expense of a reduced CFL condition. Abgrall et.al. \cite{abgrall2024new} developed a novel well-balanced active flux method for the one-dimensional shallow water equations. Furthermore, recent work by Barsukow et.al. \cite{barsukow2026equivalence} has shown that the semi-discrete active flux and discontinuous Galerkin are the same scheme for linear problems in one and several dimensions, and in a qualified sense for nonlinear one-dimensional systems. Because of the higher polynomial degree per degree of freedom for DG, the authors show that active flux is more economical than DG for a target error threshold. Nevertheless, the two schemes exhibit significant differences in their representation of the solution, implementation, and computational cost.

In this paper, we develop a structure-preserving scheme for solving the Euler equations with gravitation, utilizing the semi-discrete active flux approach. The designed scheme aims to simultaneously maintain the generalized hydrostatic equilibrium state \eqref{hydro:euler} with high-order accuracy and preserve the positivity of the density and pressure. 
For the evolution of the cell average, we adopt a special source term treatment, inspired by \cite{li2016highfve}, to preserve the well-balanced property. The gravitational potential is reformulated into an equivalent form based on the corresponding hydrostatic equilibrium solution.
For the evolution of the point values, we employ the strategy of Lax-Friedrichs flux vector splitting developed in \cite{duan2025active}, which shows better performance than the Jacobian splitting \cite{abgrall2025semi,abgrall2023extensions} for strong discontinuities. 
Drawing inspiration from \cite{xing2013highfde}, where the source term mimics the same weighted essentially non-oscillatory approximation to the flux term, we apply the same finite difference operator to both the flux term and the source term. The well-balanced property is also achieved by making a minor modification to the Lax-Friedrichs flux vector splitting. Besides the well-balancedness, the constructed scheme is also positivity-preserving. To address nonphysical solutions or spurious oscillations near discontinuities generated by the high-order scheme, we apply a posteriori MOOD paradigm introduced in \cite{clain2011high,abgrall2024new} to both the cell average and point values. 
Theoretical analysis and numerical experiments demonstrate the well-balanced and positivity-preserving nature of the proposed scheme, as well as its ability to achieve high resolution.

An outline of this paper is as follows. 
Section \ref{se:wb} introduces the well-balanced active flux method for the Euler equations with gravity, which can preserve both the isothermal and polytropic equilibrium simultaneously. In addition to the well-balanced property, the proposed scheme also preserves the positivity of the density and pressure. In Section \ref{se:nu}, numerical experiments in different circumstances are illustrated to validate all the mathematical properties of the active flux method.  In the end, some concluding remarks are provided in Section \ref{se:co}.

\section{Well-balanced active flux method}\label{se:wb}
In this section, we present the well-balanced and positivity-preserving active flux method for the Euler equations with gravitation, with the aim of preserving the general hydrostatic states \eqref{hydro:euler}. Here we restrict ourselves to two dimensions, the Euler equations \eqref{euler2D} take the form of
\begin{equation}\label{hbl}
  \boldsymbol{U}_t + \boldsymbol{F}(\boldsymbol{U})_x + \boldsymbol{G}(\boldsymbol{U})_y = \boldsymbol{S}(\boldsymbol{U},\phi),
\end{equation}
with
\begin{equation}\label{co}
  \boldsymbol{F}(\boldsymbol{U}) = \left[\begin{array}{c}
    \rho u    \\
    \rho u^2+p\\
    \rho uv \\
    (E+p)u
    \end{array} \right], \
  \boldsymbol{G}(\boldsymbol{U}) = \left[\begin{array}{c}
      \rho v    \\
      \rho uv \\
      \rho v^2+p\\ 
      (E+p)v
      \end{array} \right],\
  \boldsymbol{S}(\boldsymbol{U},\phi) = \left[\begin{array}{c}
    0    \\
    -\rho \phi_x \\
    -\rho \phi_y\\ 
    -\rho u \phi_x -\rho v\phi_y
    \end{array} \right].
\end{equation}

Assuming that the computational domain $\Omega$ is divided into $Nx\times Ny$ rectangular cells $ I_{ij}$, i.e.
\begin{equation*}
  \Omega = \bigcup  _{i=1}^{Nx}\bigcup _{j=1}^{Ny}{I_{ij}} = \bigcup  _{i=1}^{Nx}\bigcup _{j=1}^{Ny} [x_{i-\frac{1}{2}},x_{i+\frac{1}{2}}] \times [y_{j-\frac{1}{2}},y_{j+\frac{1}{2}}].
\end{equation*}
Without loss of generality, uniform meshes are considered in this work.  
The cell lengths are denoted by $\Delta x = x_{i+\frac{1}{2}}-x_{i-\frac{1}{2}}$, $\Delta y = y_{j+\frac{1}{2}}-y_{j-\frac{1}{2}}$, with corresponding cell centers located at $(x_i,y_j) = (\frac{1}{2}(x_{i-\frac{1}{2}}+x_{i+\frac{1}{2}}),\frac{1}{2}(y_{j-\frac{1}{2}}+y_{j+\frac{1}{2}}))$. The degrees of freedom involves the cell average $\overline{\boldsymbol{U}}_{ij}(t) =\frac{1}{\Delta x\Delta y}\int_{I_{ij}}\boldsymbol{U}(x,y,t) \ dxdy$ as well as the point values $\boldsymbol{U}_{\sigma}$ located at the cell interfaces. Here $\boldsymbol{U}(x,y,t)$ is the high-order numerical solution, $\sigma=(i+\frac{1}{2},j), (i,j+\frac{1}{2}),(i+\frac{1}{2},j+\frac{1}{2}) $. Specifically, the face-centered values $\boldsymbol{U}_{i+\frac{1}{2},j}(t)=\boldsymbol{U}(x_{i+\frac{1}{2}},y_j,t),\boldsymbol{U}_{i,j+\frac{1}{2}}=\boldsymbol{U}(x_i,y_{j+\frac{1}{2}},t)$ and the values at the corner $\boldsymbol{U}_{i+\frac{1}{2},j+\frac{1}{2}}=\boldsymbol{U}(x_{i+\frac{1}{2}},y_{j+\frac{1}{2}},t)$ are included for the third-order well-balanced active flux method.

\subsection{Well-balanced evolution of cell average}
Following the procedure of \cite{duan2025active}, the well-balanced evolution of the cell average can be expressed as:
\begin{equation}\label{ave}
  \dfrac{d}{dt}\overline{\boldsymbol{U}}_{ij} = -\dfrac{1}{\Delta x}(\widehat{\boldsymbol{F}}_{i+\frac{1}{2},j} - \widehat{\boldsymbol{F}}_{i-\frac{1}{2},j}) - \dfrac{1}{\Delta y}(\widehat{\boldsymbol{G}}_{i,j+\frac{1}{2}} - \widehat{\boldsymbol{G}}_{i,j-\frac{1}{2}})+\dfrac{1}{\Delta x_i \Delta y_j}\overline{\boldsymbol{S}}_{ij},
\end{equation}
with the numerical fluxes defined as follows,
\begin{equation}\label{flux}
  \begin{aligned}
  \widehat{\boldsymbol{F}}_{i+\frac{1}{2},j} &= \dfrac{1}{\Delta y}\int_{y_{j-\frac{1}{2}}}^{y_{j+\frac{1}{2}}} \boldsymbol{F}(\boldsymbol{U}(x_{i+\frac{1}{2}},y)) \ dy,\\
  \widehat{\boldsymbol{G}}_{i,j+\frac{1}{2}} &= \dfrac{1}{\Delta x}\int_{x_{i-\frac{1}{2}}}^{x_{i+\frac{1}{2}}} \boldsymbol{G}(\boldsymbol{U}(x,y_{j+\frac{1}{2}})) \ dx.
  \end{aligned}
\end{equation}
To attain third-order accuracy, Simpson's formula is employed to evaluate \eqref{flux}, yielding
\begin{equation}\label{flux_3rd}
  \widehat{\boldsymbol{F}}_{i+\frac{1}{2},j} = \dfrac{1}{6}(\boldsymbol{F}(\boldsymbol{U}_{i+\frac{1}{2},j-\frac{1}{2}}) + 4\boldsymbol{F}(\boldsymbol{U}_{i+\frac{1}{2},j}) + \boldsymbol{F}(\boldsymbol{U}_{i+\frac{1}{2},j+\frac{1}{2}})).
\end{equation}
The numerical flux $\widehat{\boldsymbol{G}}_{i,j+\frac{1}{2}}$ can be obtained in an analogous manner.

Moreover, $\overline{\boldsymbol{S}}_{ij}$ is the high-order approximation of $\int_{I_{ij}} \boldsymbol{S}(\boldsymbol{U}(x,y,t),\phi(x,y))\ dxdy$. Direct computation of this quantity generally leads to a loss of the well-balanced property of the scheme. Therefore, a special treatment of the source term is necessary to maintain the equilibrium state.  
Suppose that the target hydrostatic equilibrium $\{\rho^e(x,y),\boldsymbol{u}^e(x,y),p^e(x,y)\}$ is known a priori, which implies the following relation holds:
\begin{equation}\label{u_eq}
  \nabla p^e = \rho^e \nabla \phi, \quad \boldsymbol{u}^e=0.
\end{equation}
Taking the second component of the source term as an example, which is denoted by $ \boldsymbol{S}_{ij}^{(2)}$. We follow the similar technique in \cite{wu2021uniformly,li2016highfve} to reformulate and decompose its integral as 
\begin{equation}\label{ave_s_2}
  \begin{aligned}
    \overline{\boldsymbol{S}}_{ij}^{(2)} &= \int_{I_{ij}} -\rho \phi_x \ dx dy = \int_{I_{ij}} \dfrac{\rho}{\rho^e} p_x^e \ dx dy = \int_{I_{ij}} \left( \dfrac{\rho}{\rho^e} -  \dfrac{\bar\rho_{ij}}{\bar\rho^e_{ij}} + \dfrac{\bar\rho_{ij}}{\bar\rho^e_{ij}}\right) p_x^e \  dx dy \\
  & = \int_{I_{ij}} \left( \dfrac{\rho}{\rho^e} -  \dfrac{\bar\rho_{ij}}{\bar\rho^e_{ij}} \right) p_x^e \  dx dy + \dfrac{\bar\rho_{ij}}{\bar\rho^e_{ij}}\int_{y_{j-\frac{1}{2}}}^{y_{j+\frac{1}{2}}} ( p^e_{i+\frac{1}{2},j} - p^e_{i-\frac{1}{2},j})\ dy. 
  \end{aligned}
\end{equation}
Here $\bar{\cdot}_{ij}$ denotes the cell average of the corresponding quantity in cell $I_{ij}$. The first integral of the final identity can be approximated by the standard Gauss quadrature rule, with $\{(x_{i_{\mu}},y_{j_{\nu}}),\omega_{\mu\nu}\}$ being the Gauss quadrature nodes and associated weights. The second integral is evaluated using a Simpson formula similar to that in \eqref{flux_3rd}. For ease of understanding, we have
\begin{equation}\label{ave_ss_2}
  \begin{aligned}
  \overline{\boldsymbol{S}}_{ij}^{(2)} &\approx \sum_{\mu,\nu}\omega_{\mu\nu}
   \left( \dfrac{\rho_{i_{\mu},j_{\nu}}}{\rho^e_{i_{\mu},j_{\nu}}} -  \dfrac{\bar\rho_{ij}}{\bar\rho^e_{ij}} \right) p_x^e(x_{i_{\mu}},y_{j_{\nu}}) \\
   & + \dfrac{\bar\rho_{ij}}{\bar\rho^e_{ij}}\dfrac{\Delta y}{6}\left( (p^e_{i+\frac{1}{2},j-\frac{1}{2}}+4p^e_{i+\frac{1}{2},j}+p^e_{i+\frac{1}{2},j+\frac{1}{2}} )-  (p^e_{i-\frac{1}{2},j-\frac{1}{2}}+4p^e_{i-\frac{1}{2},j}+p^e_{i-\frac{1}{2},j+\frac{1}{2}} )\right).
  \end{aligned}
\end{equation}
The third and forth part of the source term $ \overline{\boldsymbol{S}}_{ij}^{(3)},\overline{\boldsymbol{S}}_{ij}^{(4)}$ can be treated in a similar manner.

\subsection{Well-balanced evolution of point values}
Apart from the evolution of the cell average, third-order active flux method also requires the evolution of the point values at the cell interfaces, which yields a consistent approximation of the flux derivative. In this paper, we update the point values by utilizing the flux vector splitting proposed in \cite{duan2025active}, which is given by 
\begin{align}
  \dfrac{d}{dt}\boldsymbol{U}_{i+\frac{1}{2},j+\frac{1}{2}} &= -[\mathcal{D}_1^+ \boldsymbol{F}^+ + \mathcal{D}_1^- \boldsymbol{F}^-]_{i+\frac{1}{2},j+\frac{1}{2}} - [\mathcal{D}_2^+ \boldsymbol{G}^+ + \mathcal{D}_2^- \boldsymbol{G}^-]_{i+\frac{1}{2},j+\frac{1}{2}} + \boldsymbol{S}_{i+\frac{1}{2},j+\frac{1}{2}},\label{cell_1} \\
  \dfrac{d}{dt}\boldsymbol{U}_{i+\frac{1}{2},j} &= -[\mathcal{D}_1^+ \boldsymbol{F}^+ + \mathcal{D}_1^- \boldsymbol{F}^-]_{i+\frac{1}{2},j} - [ \mathcal{D}_2 \boldsymbol{G}]_{i+\frac{1}{2},j}+ \boldsymbol{S}_{i+\frac{1}{2},j},\label{cell_2} \\
  \dfrac{d}{dt}\boldsymbol{U}_{i,j+\frac{1}{2}} &= -[\mathcal{D}_1 \boldsymbol{F}]_{i,j+\frac{1}{2}} - [\mathcal{D}_2^+ \boldsymbol{G}^+ + \mathcal{D}_2^- \boldsymbol{G}^-]_{i,j+\frac{1}{2}} + \boldsymbol{S}_{i,j+\frac{1}{2}}.\label{cell_3}
\end{align}
Here the flux $\boldsymbol{F} = \boldsymbol{F}^+ + \boldsymbol{F}^-$ and satisfies $\lambda(\frac{\partial\boldsymbol{F}^+}{\partial \boldsymbol{U}})\geq 0,\lambda(\frac{\partial\boldsymbol{F}^-}{\partial \boldsymbol{U}})\leq 0$. Considering the well-balanced property, we adopt the Lax-Friedrichs flux vector splitting with slight modification, where $\boldsymbol{F}^{\pm}$ is defined by
\begin{equation}\label{llf}
  \boldsymbol{F}^{\pm}=\dfrac{1}{2}(\boldsymbol{F}(\boldsymbol{U})\pm \alpha \widetilde{\boldsymbol{U}}), \quad \widetilde{\boldsymbol{U}} = \boldsymbol{U}-\boldsymbol{U}^e. 
\end{equation}
$\boldsymbol{U}^e$ is the conservative variables at the equilibrium state, which can be obtained from \eqref{u_eq} and  
$\alpha$ is taken as the global maximum of the spectral radius of $\frac{\partial\boldsymbol{F}}{\partial \boldsymbol{U}}$ over the computation domain. The flux $\boldsymbol{G}^{\pm}$ can be defined analogously. The cell-centered point value is obtained through the bi-parabolic reconstruction using the cell average and point values at the cell interfaces,
\begin{small}
\begin{equation}\label{middle_point}
  \boldsymbol{U}_{ij} = \dfrac{1}{16}\left[36\overline{\boldsymbol{U}}_{ij} - 4(\boldsymbol{U}_{i-\frac{1}{2},j}+\boldsymbol{U}_{i+\frac{1}{2},j}+\boldsymbol{U}_{i,j-\frac{1}{2}}+\boldsymbol{U}_{i,j+\frac{1}{2}})-(\boldsymbol{U}_{i-\frac{1}{2},j-\frac{1}{2}}+\boldsymbol{U}_{i-\frac{1}{2},j+\frac{1}{2}}+\boldsymbol{U}_{i+\frac{1}{2},j-\frac{1}{2}}+\boldsymbol{U}_{i+\frac{1}{2},j+\frac{1}{2}})\right]
\end{equation}
\end{small}
More details can be found in \cite{abgrall2025semi}. We can easily obtain the fluxes on the edges and the cell-centered point value with the definition of \eqref{llf}.
Furthermore, $\mathcal D_1^{\pm}, \mathcal D_2^{\pm}$ are the finite difference operators acquired from differentiating the bi-parabolic interpolation of $\boldsymbol{F}^{\pm},\boldsymbol{G}^{\pm}$ in the cell $I_{ij}$, 
\begin{equation}\label{operator}
  \begin{aligned}
    [\mathcal{D}_1^+ \boldsymbol{F}^+]_{i+\frac{1}{2},j_{\beta}} &= \dfrac{1}{\Delta x}(\boldsymbol{F}^{+}_{i-\frac{1}{2},j_{\beta}} - 4\boldsymbol{F}^{+}_{i,j_{\beta}}+3\boldsymbol{F}^{+}_{i+\frac{1}{2},j_{\beta}}),\\
    [\mathcal{D}_1^- \boldsymbol{F}^-]_{i+\frac{1}{2},j_{\beta}} &= \dfrac{1}{\Delta x}(-3\boldsymbol{F}^{-}_{i+\frac{1}{2},j_{\beta}}+4\boldsymbol{F}^{-}_{i+1,j_{\beta}}-\boldsymbol{F}^{-}_{i+\frac{3}{2},j_{\beta}}),\\
    [\mathcal{D}_1\boldsymbol{F}]_{i,j_{\beta}}& = \dfrac{1}{\Delta x}(\boldsymbol{F}_{i+\frac{1}{2},j_{\beta}} - \boldsymbol{F}_{i-\frac{1}{2},j_{\beta}}), \\
    [\mathcal{D}_2^+ \boldsymbol{G}^+]_{i_{\beta},j+\frac{1}{2}} &=\dfrac{1}{\Delta y}(\boldsymbol{G}^{+}_{i_{\beta},j-\frac{1}{2}} - 4\boldsymbol{G}^{+}_{i_{\beta},j}+3\boldsymbol{G}^{+}_{i_{\beta},j+\frac{1}{2}}),\\
    [\mathcal{D}_2^- \boldsymbol{G}^-]_{i_{\beta},j+\frac{1}{2}} &=\dfrac{1}{\Delta y}(-3\boldsymbol{G}^{-}_{i_{\beta},j+\frac{1}{2}}+4\boldsymbol{G}^{-}_{i_{\beta},j+1}-\boldsymbol{G}^{-}_{i_{\beta},j+\frac{3}{2}}),\\
    [\mathcal{D}_2 \boldsymbol{G}]_{i_{\beta},j}  &=\dfrac{1}{\Delta y}(\boldsymbol{G}_{i_{\beta},j+\frac{1}{2}} - \boldsymbol{G}_{i_{\beta},j-\frac{1}{2}}),
  \end{aligned}
\end{equation}
where $i_{\beta} = \{i-\frac{1}{2},i,i+\frac{1}{2}\}$, $j_{\beta} = \{j-\frac{1}{2},j,j+\frac{1}{2}\}$. 

For the source term approximation, we employ the relation \eqref{u_eq} and reformulate the second component of the source term as
\begin{equation}\label{cell_s_2} 
    \boldsymbol{S}_{i+\frac{1}{2},j+\frac{1}{2}}^{(2)} = [-\rho \phi_x]_{i+\frac{1}{2},j+\frac{1}{2}} = \left(\frac{\rho}{\rho^e}\right)_{i+\frac{1}{2},j+\frac{1}{2}}(p^e_x)_{i+\frac{1}{2},j+\frac{1}{2}}.
\end{equation}
Inspired by \cite{xing2013highfde} and in order to preserve the well-balanced property, we split the derivatives in the source term as $\nabla p^e = \nabla(\frac{1}{2}p^e) +  \nabla(\frac{1}{2}p^e)$, and apply the same finite difference operator to approximate them. We have
\begin{equation}\label{cell_s}
  \begin{aligned}
    \boldsymbol{S}_{i+\frac{1}{2},j+\frac{1}{2}}^{(2)} &=\left(\frac{\rho}{\rho^e}\right)_{i+\frac{1}{2},j+\frac{1}{2}}\left[\mathcal{D}^+_1\left(\frac{1}{2}p^e\right)+\mathcal{D}^-_1\left(\frac{1}{2}p^e\right)\right]_{i+\frac{1}{2},j+\frac{1}{2}},\\
    \boldsymbol{S}_{i+\frac{1}{2},j+\frac{1}{2}}^{(3)} &=\left(\frac{\rho}{\rho^e}\right)_{i+\frac{1}{2},j+\frac{1}{2}}\left[\mathcal{D}^+_2\left(\frac{1}{2}p^e\right)+\mathcal{D}^-_2\left(\frac{1}{2}p^e\right)\right]_{i+\frac{1}{2},j+\frac{1}{2}},\\
    \boldsymbol{S}_{i+\frac{1}{2},j+\frac{1}{2}}^{(4)} &=\left(\frac{\rho u}{\rho^e}\right)_{i+\frac{1}{2},j+\frac{1}{2}}\left[\mathcal{D}^+_1\left(\frac{1}{2}p^e\right)+\mathcal{D}^-_1\left(\frac{1}{2}p^e\right)\right]_{i+\frac{1}{2},j+\frac{1}{2}} \\
    &+ \left(\frac{\rho v}{\rho^e}\right)_{i+\frac{1}{2},j+\frac{1}{2}}\left[\mathcal{D}^+_2\left(\frac{1}{2}p^e\right)+\mathcal{D}^-_2\left(\frac{1}{2}p^e\right)\right]_{i+\frac{1}{2},j+\frac{1}{2}}.\\
  \end{aligned}
\end{equation}
Such choice of splitting corresponds to the Lax-Friedrichs flux vector splitting and plays a crucial role in a well-balanced manner. The face-centered source term  in \eqref{cell_2} is written as
\begin{equation}\label{cell_sl}
  \begin{aligned}
    \boldsymbol{S}_{i+\frac{1}{2},j}^{(2)} &=\left(\frac{\rho}{\rho^e}\right)_{i+\frac{1}{2},j}\left[\mathcal{D}^+_1\left(\frac{1}{2}p^e\right)+\mathcal{D}^-_1\left(\frac{1}{2}p^e\right)\right]_{i+\frac{1}{2},j},\\
    \boldsymbol{S}_{i+\frac{1}{2},j}^{(3)} &=\left(\frac{\rho}{\rho^e}\right)_{i+\frac{1}{2},j}\left[\mathcal{D}_2\left(p^e\right)\right]_{i+\frac{1}{2},j},\\
    \boldsymbol{S}_{i+\frac{1}{2},j}^{(4)} &=\left(\frac{\rho u}{\rho^e}\right)_{i+\frac{1}{2},j}\left[\mathcal{D}^+_1\left(\frac{1}{2}p^e\right)+\mathcal{D}^-_1\left(\frac{1}{2}p^e\right)\right]_{i+\frac{1}{2},j}+\left(\frac{\rho v}{\rho^e}\right)_{i+\frac{1}{2},j}\left[\mathcal{D}_2\left(p^e\right)\right]_{i+\frac{1}{2},j},
  \end{aligned}
\end{equation}
and $ \boldsymbol{S}_{i,j+\frac{1}{2}}$ in \eqref{cell_3} can be obtained similarly.

\subsection{Time discretization}
For hyperbolic conservation laws,  the semi-discrete scheme is usually paired with high-order strong stability preserving Runge-Kutta (SSP-RK) or multistep time discretization methods \cite{gottlieb2001strong,shu1988total}. Here, the third-order SSP-RK method is taken for our time discretization,
\begin{equation}\label{RK3}
  \begin{aligned}
  \boldsymbol{U}^{(1)} &= \boldsymbol{U}^{n}+\Delta t \text{RHS}(\boldsymbol{U}^{n}),\\
  \boldsymbol{U}^{(2)} &= \frac{3}{4}\boldsymbol{U}^{n}+ \frac{1}{4}(\boldsymbol{U}^{(1)}+\Delta t \text{RHS}(\boldsymbol{U}^{(1)})),\\
  \boldsymbol{U}^{n+1} &= \frac{1}{3}\boldsymbol{U}^{n}+ \frac{2}{3}(\boldsymbol{U}^{(2)}+\Delta t \text{RHS}(\boldsymbol{U}^{(2)})),
  \end{aligned}
\end{equation}
where $\text{RHS}$ are the spatial operators given on the right hand side of the schemes \eqref{ave}, \eqref{cell_1}-\eqref{cell_3}.

\subsection{Well-balanced property}
We present the well-balanced property for the semi-discrete and fully-discrete active flux method in the following theorems.

\begin{theorem}\label{the:semi}
  The semi-discrete active flux method given by  \eqref{ave}, \eqref{cell_1}-\eqref{cell_3} is well-balanced for the hydrostatic equilibrium state \eqref{hydro:euler}.
\end{theorem}

\begin{proof}
  Assuming the initial conditions correspond to an exact hydrostatic equilibrium state, which means 
  \begin{equation}\label{steady_ini}
    \rho(x,y) = \rho^e(x,y), \ \boldsymbol{u}(x,y) = 0, \ p(x,y) = p^e(x,y), \ E(x,y) = \frac{1}{\gamma-1}p^e(x,y).
  \end{equation}
  Since $\boldsymbol{u}=0$, the fluxes term reduce to
  \begin{equation}\label{flux_eq}
    \boldsymbol{F}(\boldsymbol{U}) = (0,p,0,0)^T, \quad \boldsymbol{G}(\boldsymbol{U}) = (0,0,p,0)^T.
  \end{equation}
  It is straightforward to verify that for the cell average \eqref{ave}, the well-balanced property holds for the first and fourth equations, since both the flux and the source term approximations become zero in equilibrium.
  From the steady state  \eqref{steady_ini}, we have $\bar \rho_{ij} = \bar\rho_{ij}^e$. Therefore, the first term of \eqref{ave_ss_2} vanishes in equilibrium, and
  \begin{equation}\label{ave_s2_eq}
    \overline{\boldsymbol{S}}_{ij}^{(2)} = \dfrac{\Delta y}{6}\left( (p^e_{i+\frac{1}{2},j-\frac{1}{2}}+4p^e_{i+\frac{1}{2},j}+p^e_{i+\frac{1}{2},j+\frac{1}{2}} )-  (p^e_{i-\frac{1}{2},j-\frac{1}{2}}+4p^e_{i-\frac{1}{2},j}+p^e_{i-\frac{1}{2},j+\frac{1}{2}} )\right).
  \end{equation}
  By substituting \eqref{flux_eq} and \eqref{ave_s2_eq} into \eqref{ave}, we get
  \begin{align*}
    \dfrac{d}{dt}\overline{\boldsymbol{U}}_{ij}^{(2)} &= -\dfrac{1}{\Delta x}(\widehat{\boldsymbol{F}}_{i+\frac{1}{2},j}^{(2)} - \widehat{\boldsymbol{F}}_{i-\frac{1}{2},j}^{(2)}) - \dfrac{1}{\Delta y}(\widehat{\boldsymbol{G}}_{i,j+\frac{1}{2}}^{(2)} - \widehat{\boldsymbol{G}}_{i,j-\frac{1}{2}}^{(2)}) + \dfrac{1}{\Delta x \Delta y} \overline{\boldsymbol{S}}_{ij}^{(2)} \\
    &=-\dfrac{1}{\Delta x}\dfrac{1}{6}\left( (p_{i+\frac{1}{2},j-\frac{1}{2}}+4p_{i+\frac{1}{2},j}+p_{i+\frac{1}{2},j+\frac{1}{2}} )-  (p_{i-\frac{1}{2},j-\frac{1}{2}}+4p_{i-\frac{1}{2},j}+p_{i-\frac{1}{2},j+\frac{1}{2}} )\right) \\
    &+\dfrac{1}{\Delta x \Delta y}\dfrac{\Delta y}{6}\left( (p^e_{i+\frac{1}{2},j-\frac{1}{2}}+4p^e_{i+\frac{1}{2},j}+p^e_{i+\frac{1}{2},j+\frac{1}{2}} )-  (p^e_{i-\frac{1}{2},j-\frac{1}{2}}+4p^e_{i-\frac{1}{2},j}+p^e_{i-\frac{1}{2},j+\frac{1}{2}} )\right) \\
    &= 0.
  \end{align*}
  The last equality follows from the equilibrium state \eqref{steady_ini}. Likewise, we can obtain $\frac{d}{dt}\overline{\boldsymbol{U}}_{ij}^{(3)}=0$. Then $\frac{d}{dt}\overline{\boldsymbol{U}}_{ij}=0$ follows.

  For the point values, thanks to the steady state  \eqref{steady_ini}, we have $\widetilde{\boldsymbol{U}}=0$, thus the fluxes become 
  $$\boldsymbol{F}^{\pm} = \dfrac{1}{2}\boldsymbol{F}(\boldsymbol{U})=(0,\frac{1}{2}p,0,0)^T,\quad \boldsymbol{G}^{\pm} = \dfrac{1}{2}\boldsymbol{G}(\boldsymbol{U})=(0,0,\frac{1}{2}p,0)^T.$$ 
  Analogously, it is trivial that the well-balanced property holds for the mass and energy equations, we mainly focus on the momentum equation. 
  Combined with the fact that $(\rho/\rho^e)_{i+\frac{1}{2},j+\frac{1}{2}}=1$ and the same finite difference operator performed on the source term, we have 
  \begin{align*}
    \dfrac{d}{dt}\boldsymbol{U}_{i+\frac{1}{2},j+\frac{1}{2}}^{(2)} &= -[\mathcal{D}_1^+ \boldsymbol{F}^{(2),+} + \mathcal{D}_1^- \boldsymbol{F}^{(2),-}]_{i+\frac{1}{2},j+\frac{1}{2}} - [\mathcal{D}_2^+ \boldsymbol{G}^{(2),+} + \mathcal{D}_2^- \boldsymbol{G}^{(2),-}]_{i+\frac{1}{2},j+\frac{1}{2}} + \boldsymbol{S}_{i+\frac{1}{2},j+\frac{1}{2}}^{(2)}\\
    &=-\left[\mathcal{D}^+_1\left(\frac{1}{2}p\right)+\mathcal{D}^-_1\left(\frac{1}{2}p\right)\right]_{i+\frac{1}{2},j+\frac{1}{2}}+\left[\mathcal{D}^+_1\left(\frac{1}{2}p^e\right)+\mathcal{D}^-_1\left(\frac{1}{2}p^e\right)\right]_{i+\frac{1}{2},j+\frac{1}{2}}\\
    &=0.
  \end{align*}
  The third component, as well as the remaining face-centered point values \eqref{cell_2} and \eqref{cell_3}, can be acquired in a similar manner. This completes the proof of the well-balanced property of our scheme.
\end{proof}

A straightforward induction yields the well-balancedness of the full-discrete scheme, and we omit the detailed proof in this context.
\begin{theorem}\label{the:fully}
  The fully-discrete active flux method, given by \eqref{ave}, \eqref{cell_1}-\eqref{cell_3} coupled with third-order SSP-RK time discretization, is well-balanced for the hydrostatic equilibrium state \eqref{hydro:euler}.
\end{theorem}

\subsection{Positivity-preserving property} 
In physics, the conservative variables should belong to the set of admissible state, which is defined by
\begin{equation}\label{admis}
  \mathcal{G} := \left\{ \boldsymbol{U} = (\rho,\rho\boldsymbol{u},E)^T |\ \rho>0, \ p =p(\boldsymbol{U})= (\gamma-1)\left(E-\frac{||\rho\boldsymbol{u}||^2}{2\rho}\right)>0\right\}.
\end{equation}
This property, however, is not automatically guaranteed by numerical schemes, which may produce negative density or pressure if we do not take special treatment in low density and low pressure regions. It is also known that high-order schemes are prone to numerical oscillations when containing discontinuous solutions. Such numerical instability might cause the breakdown of numerical simulations. Therefore, in addition to achieving well-balanced property, it is of vital importance to preserve the positivity of the density and pressure while effectively controlling spurious oscillations.

For the proposed active flux method, we should keep both the cell average $\overline{\boldsymbol{U}}_{ij}\in \mathcal{G}$ and the point values $\boldsymbol{U}_{\sigma}\in \mathcal{G}$. To this end, we apply a simplified version of the MOOD paradigm which was used in \cite{clain2011high,diot2012improved,abgrall2024new,abgrall2023combination}. 
The MOOD procedure is divided into three parts. First, the candidate solution is computed by the well-balanced third-order active flux method \eqref{ave}, \eqref{cell_1}-\eqref{cell_3}. Then we check whether the recomputation is needed in each cell according to a set of detection criteria. Finally, if a cell is flagged as a 'trouble' cell, we use a more robust lower-order scheme for this specific cell.

\subsubsection{Detection procedure}
Let $(\overline{\boldsymbol{U}}_{ij}^n, \boldsymbol{U}_{\sigma}^n)$ be the solution at time $t_n$. After each inner stage of Runge-Kutta method, the candidate solution is denoted as $(\overline{\boldsymbol{U}}_{ij}^{*}, \boldsymbol{U}_{\sigma}^{*})$. The following criteria are applied to detect trouble cells for the cell average,
\begin{enumerate}
  \setlength{\parskip}{0.0em}
  \item Computer admissible detector: We check if all the variables are numbers, i.e. not NaN or Inf. If not, we flag the element .FALSE. and proceed to the next element. Otherwise, we set the element .TRUE..
  \item Physically admissible detector: We check if $\overline{\boldsymbol{U}}_{ij}^{*}\in \mathcal{G}$. If this holds, we set this element .TRUE. and proceed to check the next criteria, otherwise the element is flagged .FALSE.. 
  \item Locally constant detection: We check whether the solution is not locally constant at time $t_n$. To be specific, we compare the values of the degrees of freedom in the element $I_{ij}$ with those in the neighboring elements sharing a face. A procedure similar to that in \cite{abgrall2024new} is followed to avoid detecting a fake maximum principle. If the solution is locally constant, no modification of the element is required and the corresponding flag is set to .TRUE..
  \item Numerical admissible detection: We check whether a new extremum is produced or not via the discrete maximum principle. Here we test for the density, and possibly the pressure. Let $\xi$ be the function on which we perform the test, and $\mathcal{V}_{ij}$ be the set containing $I_{ij}$ together with its neighboring cells that share a face. Following the procedure in \cite{dumbser2014posteriori}, if
  \begin{equation*}
    \bar\xi_{ij}^* \in[\min_{l\in\mathcal{V}_{ij}}\xi_l^n - \delta, \max_{l\in\mathcal{V}_{ij}}\xi_l^n+\delta], \quad \delta =\max(\delta_0,\epsilon(\max_{l\in\mathcal{V}_{ij}}\xi_l^n - \min_{l\in\mathcal{V}_{ij}}\xi_l^n)),
  \end{equation*}
  with $\delta_0=10^{-4}, \epsilon = 10^{-3}$, we set the element .TRUE.. Otherwise, we employ the technique introduced in \cite{vilar2019posteriori} to avoid misidentifying smooth extrema as spurious ones, thereby preserving the accuracy of the scheme.
\end{enumerate}

The similar detection procedure is applied to the point values as well. 

\subsubsection{Recomputation procedure}
For any cell in which the candidate numerical high-order solution fails to respect the constraints, we switch to a parachute first-order scheme. Notably, not all first-order schemes possess both the well-balanced and positivity-preserving properties simultaneously. In order not to destroy these properties, the parachute scheme used here is a first-order scheme specially tailored for the Euler equations with gravity. 

For the cell average, we employ the piecewise constant structure-preserving DG scheme proposed in \cite{wu2021uniformly}, where the test function is taken to be one. More details can be found in \cite{wu2021uniformly}.
Analogous treatments are applied to the point values, which read as
\begin{align} 
  \dfrac{d}{dt}\boldsymbol{U}_{i+\frac{1}{2},j} &= - \frac{1}{\Delta x}(\widehat{\boldsymbol{F}}^L_{i+1,j} - \widehat{\boldsymbol{F}}^L_{i,j}) - \frac{1}{\Delta y}(\widehat{\boldsymbol{G}}^L_{i+\frac{1}{2},j+\frac{1}{2}} - \widehat{\boldsymbol{G}}^L_{i+\frac{1}{2},j-\frac{1}{2}}) + \boldsymbol{S}^L_{i+\frac{1}{2},j},\label{1st_x}\\
  \dfrac{d}{dt}\boldsymbol{U}_{i,j+\frac{1}{2}} &= -\frac{1}{\Delta x}(\widehat{\boldsymbol{F}}^L_{i+\frac{1}{2},j+\frac{1}{2}} - \widehat{\boldsymbol{F}}^L_{i-\frac{1}{2},j+\frac{1}{2}}) - \frac{1}{\Delta y}(\widehat{\boldsymbol{G}}^L_{i,j+1} - \widehat{\boldsymbol{G}}^L_{i,j}) + \boldsymbol{S}^L_{i,j+\frac{1}{2}},\label{1st_y}\\
  \dfrac{d}{dt}\boldsymbol{U}_{i+\frac{1}{2},j+\frac{1}{2}} &= -\frac{1}{\Delta x}(\widehat{\boldsymbol{F}}^L_{i+1,j+\frac{1}{2}} - \widehat{\boldsymbol{F}}^L_{i,j+\frac{1}{2}}) - \frac{1}{\Delta y}(\widehat{\boldsymbol{G}}^L_{i+\frac{1}{2},j+1} - \widehat{\boldsymbol{G}}^L_{i+\frac{1}{2},j}) + \boldsymbol{S}^L_{i+\frac{1}{2},j+\frac{1}{2}}. \label{1st_c}
\end{align} 
Here the superscript ${}^L$ means the lower-order numerical fluxes or source term. 
\begin{equation}\label{flux_1st_p}
  \begin{aligned}
  \widehat{\boldsymbol{F}}^L_{i+1,j} = \boldsymbol{F}^{hllc}\left(\frac{p^{e,*}_{i+1,j}}{ p^e_{i+\frac{1}{2},j}}{\boldsymbol{U}}_{i+\frac{1}{2},j}, \frac{p^{e,*}_{i+1,j}}{ p^e_{i+\frac{3}{2},j}}{\boldsymbol{U}}_{i+\frac{3}{2},j} \right), 
  \widehat{\boldsymbol{G}}^L_{i+\frac{1}{2},j+\frac{1}{2}} = \boldsymbol{G}^{hllc}\left(\frac{p^{e,*}_{i+\frac{1}{2},j+\frac{1}{2}}}{p^e_{i+\frac{1}{2},j}}\boldsymbol{U}_{i+\frac{1}{2},j}, \frac{p^{e,*}_{i+\frac{1}{2},j+\frac{1}{2}}}{p^e_{i+\frac{1}{2},j+1}}\boldsymbol{U}_{i+\frac{1}{2},j+1} \right),
\end{aligned}
\end{equation}
are the modified HLLC numerical fluxes \cite{wu2021uniformly} and 
\begin{equation}\label{pe_1st}
  p^{e,*}_{i+1,j} = \frac{1}{2}( p^e_{i+\frac{1}{2},j}+  p^e_{i+\frac{3}{2},j}), \quad p^{e,*}_{i+\frac{1}{2},j+\frac{1}{2}} = \frac{1}{2}(p^e_{i+\frac{1}{2},j}+ p^e_{i+\frac{1}{2},j+1}).
\end{equation}
The source term is reformulated in the same manner as in \eqref{cell_s_2}, and then the first derivative of $p^e$ is approximated by the central difference scheme, which means
\begin{equation}\label{source_1st_p}
  \begin{aligned}
  \boldsymbol{S}^{L,(2)}_{i+\frac{1}{2},j} &=\left(\frac{\rho}{\rho^e}\right)_{i+\frac{1}{2},j}(p^e_x)_{i+\frac{1}{2},j} = \left(\frac{\rho}{\rho^e}\right)_{i+\frac{1}{2},j}\frac{p^{e}_{i+\frac{3}{2},j} - p^{e}_{i-\frac{1}{2},j}}{2\Delta x},\\
  \boldsymbol{S}^{L,(3)}_{i+\frac{1}{2},j} &=\left(\frac{\rho}{\rho^e}\right)_{i+\frac{1}{2},j}(p^e_y)_{i+\frac{1}{2},j} = \left(\frac{\rho}{\rho^e}\right)_{i+\frac{1}{2},j}\frac{p^e_{i+\frac{1}{2},j+1} - p^e_{i+\frac{1}{2},j-1}}{2\Delta y}.
  \end{aligned}
\end{equation}
The numerical fluxes and source terms for the face-centered points in $y$-direction and those at the corner points can be defined analogously. 

The well-balanced and positivity-preserving properties of the cell averages for the first-order scheme have already established in \cite{wu2021uniformly}. Some useful properties for the admissible state and HLLC numerical flux are presented in Appendix \ref{a1}. For completeness, the corresponding results for the point values are summarized in the following two theorems.
\begin{theorem}\label{wb_1st_p}
  The semi-discrete first-order scheme given by \eqref{1st_x}-\eqref{1st_c} for the point values is well-balanced for the hydrostatic equilibrium state \eqref{hydro:euler}.
\end{theorem}
\begin{proof}
  The proof is shown in Appendix \ref{a2}.
\end{proof}

\begin{theorem}\label{pp_1st_p}
  Let $\boldsymbol{U}_{\sigma}^n$ be the point values of the first-order scheme \eqref{1st_x}-\eqref{1st_c} at time $t_n$. If $\boldsymbol{U}_{\sigma}^n\in\mathcal{G}$ for all $i,j$, then we have $\boldsymbol{U}_{\sigma}^{n+1}\in\mathcal{G}$ under the CFL-type condition $\Delta t \widehat\alpha_{\sigma}\leq 1$, $\widehat\alpha_{\sigma} = \max\{\widehat\alpha_{i+\frac{1}{2},j},\widehat\alpha_{i,j+\frac{1}{2}},\widehat\alpha_{i+\frac{1}{2},j+\frac{1}{2}}\}$. Here 
  \begin{equation*}
    \widehat\alpha_{i+\frac{1}{2},j} = 2\frac{\widehat\alpha_F}{p^{e}_{i+\frac{1}{2},j}}\left(\frac{p^{e,*}_{i+1,j}+p^{e,*}_{i,j}}{\Delta x}+\frac{p^{e,*}_{i+\frac{1}{2},j+\frac{1}{2}}+p^{e,*}_{i+\frac{1}{2},j-\frac{1}{2}}}{\Delta y}\right)+\widehat\alpha_S,
  \end{equation*}
  with
  \begin{equation*}
    \widehat\alpha_F = \max\nolimits_{\boldsymbol{U}_{i\pm\frac{1}{2},j},\boldsymbol{U}_{i+\frac{3}{2},j},\boldsymbol{U}_{i+\frac{1}{2},j\pm1}}\varrho(\boldsymbol{U}), \widehat\alpha_S = \frac{1}{\rho^e_{i+\frac{1}{2},j}\sqrt{2e_{i+\frac{1}{2},j}}}\left(\big|\frac{p^{e,*}_{i+1,j}-p^{e,*}_{i,j}}{\Delta x}\big|+\big|\frac{p^{e,*}_{i+\frac{1}{2},j+\frac{1}{2}}-p^{e,*}_{i+\frac{1}{2},j-\frac{1}{2}}}{\Delta y}\big|\right).
  \end{equation*}
  $\varrho$ is the spectral radius of $\frac{\partial{\boldsymbol{F}}}{\partial \boldsymbol{U}}$, $\frac{\partial{\boldsymbol{G}}}{\partial \boldsymbol{U}}$. The other two variables $\widehat\alpha_{i,j+\frac{1}{2}},\widehat\alpha_{i+\frac{1}{2},j+\frac{1}{2}}$ can be defined analogously.
\end{theorem}
\begin{proof}
  The proof is shown in Appendix \ref{a2}.
\end{proof}
\begin{remark} 
  The cell-centered point value \eqref{middle_point}, which is used to update the point values, is not obtained as a convex combination of admissible states. Consequently, it is possible that $\boldsymbol{U}_{ij}\notin\mathcal{G}$. To address this issue, a scaling limiter \cite{zhang2011positivity} should be applied to $\boldsymbol{U}_{ij}$ at the beginning of each inner stage of SSP-RK method. Specifically, we first replace $\rho_{ij}$ by  
  \begin{equation*}
    \widehat\rho_{ij} = \theta_{ij}^1(\rho_{ij} - \overline{\rho}_{ij}) + \bar{\rho}_{ij}, \quad \theta^1_{ij} = \min\left\{1,\frac{\bar{\rho}_{ij} - \epsilon_1}{\bar{\rho}_{ij} - \rho_{ij}}\right\}, \quad \epsilon_1 = \min\{10^{-11},\bar{\rho}_{ij},\rho_{\sigma}\}.
  \end{equation*}
 We denote $\widehat{\boldsymbol{U}}_{ij} = (\widehat\rho_{ij}, \rho\boldsymbol{u}_{ij},E_{ij})^T$, the final limiting procedure is defined as follows
  \begin{equation*}
    \boldsymbol{U}_{ij}^{\text{new}} = \theta_{ij}^2(\widehat{\boldsymbol{U}}_{ij} - \overline{\boldsymbol{U}}_{ij}) + \overline{\boldsymbol{U}}_{ij},\quad \theta^2_{ij} =  \min\left\{1,\frac{p(\overline{\boldsymbol{U}}_{ij})-\epsilon_2}{p(\overline{\boldsymbol{U}}_{ij}) - p(\widehat{\boldsymbol{U}}_{ij})} \right\}, 
  \end{equation*}
  with $ \epsilon_2 = \min\{10^{-11},p(\overline{\boldsymbol{U}}_{ij}),p(\boldsymbol{U}_\sigma)\}$.
\end{remark}

\section{Numerical examples}\label{se:nu}
This section presents the implementation of our proposed structure-preserving active flux method for one- and two-dimensional Euler equations with gravitation. 
Unless otherwise specified, the CFL number is taken as 0.3 for one-dimensional case, and 0.2 for two-dimensional case. In two dimensions, except for the cell average, only the results at the corner points are reported, as those at the face-centered points exhibit similar behavior and are therefore omitted for brevity. 
We also compare our performance with that of the non-well-balanced scheme, which employs a straightforward discretization for the source term and the original Lax-Friedrichs flux vector splitting. 

\subsection{One-dimensional tests}


\begin{example}{\bf Test for one-dimensional polytropic equilibrium}\label{euler:ise_1d}
\end{example}
To demonstrate the well-balanced property of the proposed structure-preserving active flux method for the polytropic equilibrium, we consider the benchmark example introduced in \cite{LeVeque1998wave}. The corresponding steady state solution is given by
\begin{equation}\label{euler:iseini_1d}
\begin{aligned}
  \rho(x) =  \left( \rho_0^{\gamma-1}-\frac{1}{K_0}\frac{\gamma-1}{\gamma} \phi(x)\right)^{\frac{1}{\gamma-1}}, \ u(x)=0, \ p(x) = K_0\rho(x)^\gamma,
\end{aligned}
\end{equation}
where $\phi(x) = x$ is the linear gravitational potential. The ratio of specific heats is taken as $\gamma = 5/3$. The simulation is performed up to the final time $t=2$. The $L^1$ and $L^{\infty}$ errors of both the cell averages and the point values with different mesh sizes, computed in double precision, are presented in Table \ref{table:isewb_1d}. These results confirm that our proposed method is well balanced, preserving the polytropic steady state up to machine precision.

\begin{table}[htb]
  \centering
  \caption{Example \ref{euler:ise_1d}: $L^1$ and $L^{\infty}$ errors for one-dimensional polytropic hydrostatic equilibrium flow \eqref{euler:iseini_1d}.}\label{table:isewb_1d}
  \begin{tabular}{c c c c c c c c c c c c c}
   \toprule
   \multirow{2}{*}&\multirow{2}{*} {$N$}&\multicolumn{3}{c}{$L^1$ error}&\multicolumn{3}{c}{$L^{\infty}$ error}\\
    \cmidrule(lr){3-5} \cmidrule(lr){6-8}
    &&{$\rho$}&{$\rho u$}&{$E$}&{$\rho$}&{$\rho u$}&{$E$}\\
    \midrule
    \multirow{3}{*}{\text{average}}&50  & 3.73E-16 & 3.60E-16 & 3.80E-16 & 7.77E-16 & 4.06E-16 & 6.66E-16 \\
    &100 & 8.37E-16 & 3.52E-16 & 5.99E-16 & 3.33E-15 & 6.36E-16 & 1.55E-15 \\
    &200 & 1.75E-15 & 4.12E-16 & 6.35E-16 & 9.88E-15 & 8.95E-16 & 1.55E-15 \\
    \midrule
    \multirow{3}{*}{\text{point}}& 50  & 4.01E-16 & 3.60E-16 & 3.77E-16 & 9.99E-16 & 4.21E-16 & 4.44E-16 \\
    & 100 & 7.38E-16 & 3.70E-16 & 5.71E-16 & 3.00E-15 & 6.68E-16 & 1.44E-15 \\
    & 200 & 1.77E-15 & 4.24E-16 & 6.11E-16 & 9.33E-15 & 8.55E-16 & 1.33E-15 \\
    \bottomrule
  \end{tabular}
\end{table}

\begin{figure}[htb!]
  \centering
  \subfigure[pressure perturbation]{
  \centering
  \includegraphics[width=5.5cm,scale=1]{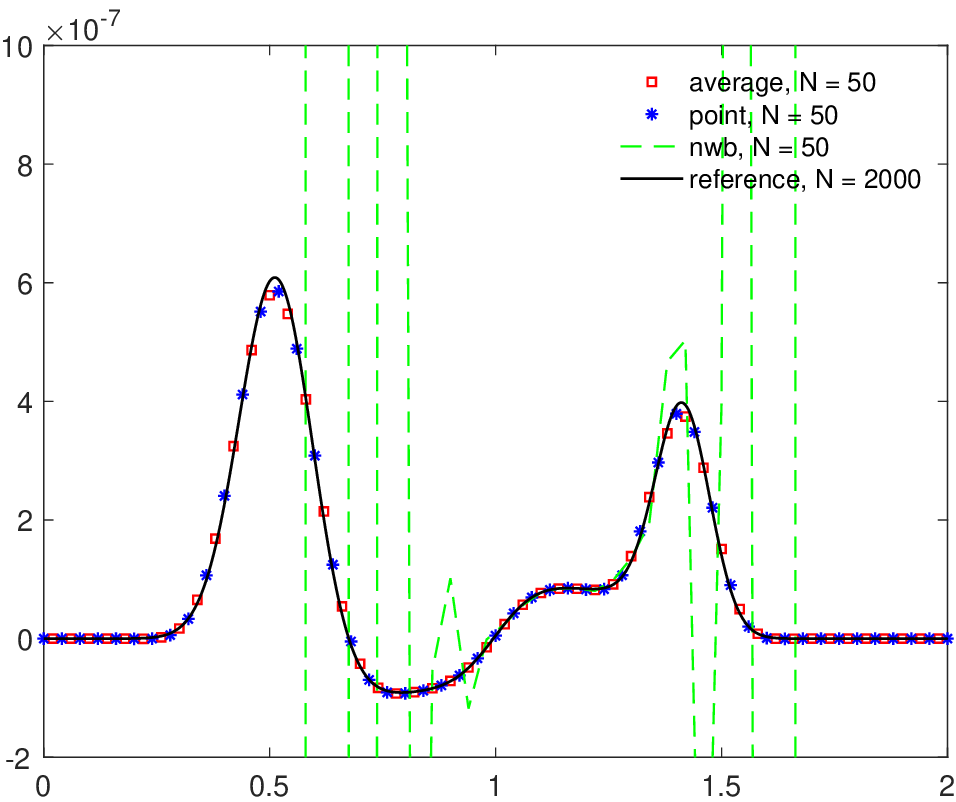}
  }\quad
  \subfigure[velocity ]{
  \centering
  \includegraphics[width=5.5cm,scale=1]{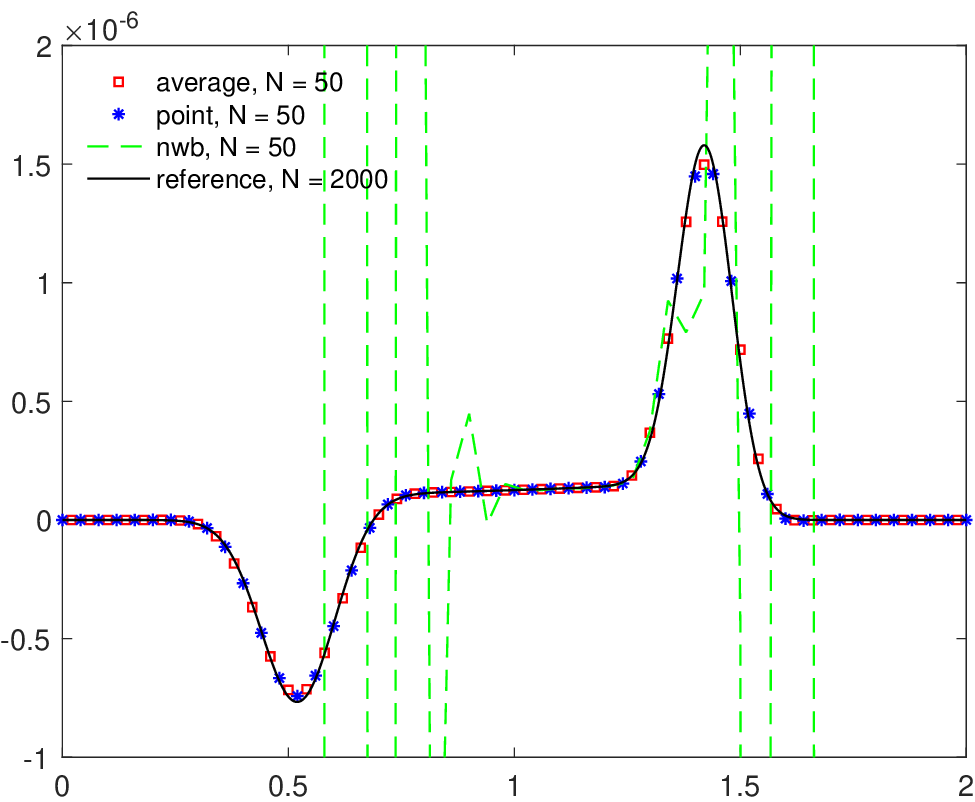}
  }
  \caption{Example \ref{euler:iseper_1d}: Pressure perturbation and velocity of the polytropic equilibrium for the small amplitude wave propagation $A=10^{-6}$. }\label{fig:persmall_p}
\end{figure}

\begin{example}{\bf Small perturbation of the polytropic equilibrium}\label{euler:iseper_1d}
\end{example}
In this example, we evaluate the ability of our proposed method to accurately propagate small perturbations near the hydrostatic equilibrium, a property which cannot be realized by a non-well-balanced scheme. Almost the same setup as in Example \ref{euler:ise_1d}, we add a small disturbance to the pressure, 
\begin{equation}\label{euler:adi_peradd}
  (\rho,u,p)(x,0) = (\rho_{eq}(x), u_{eq}(x), p_{eq}(x)+A\exp(-100(x-1)^2)).
\end{equation}
Here $(\rho_{eq}, u_{eq}, p_{eq})$ are the equilibrium functions given in \eqref{euler:iseini_1d}. We first consider a small amplitude wave propagation $A=10^{-6}$ and take the stopping time $t=0.45$. The  numerical results for the pressure perturbation and velocity computed on a coarse mesh consisting of 50 uniform cells are displayed in
Fig. \ref{fig:persmall_p}. For comparison, the same test is also performed using a non-well-balanced scheme. All results are compared with the reference solutions computed by third-order well-balanced DG scheme \cite{li2018welldge} with a much refined $2000$ mesh. As expected, the non-well-balanced scheme fails to resolve the solution accurately, while our proposed structure-preserving active flux method maintains excellent agreement with the reference solutions on a relatively coarse mesh. 

\begin{figure}[htb!]
  \centering
  \subfigure[pressure perturbation]{
  \centering
  \includegraphics[width=5.5cm,scale=1]{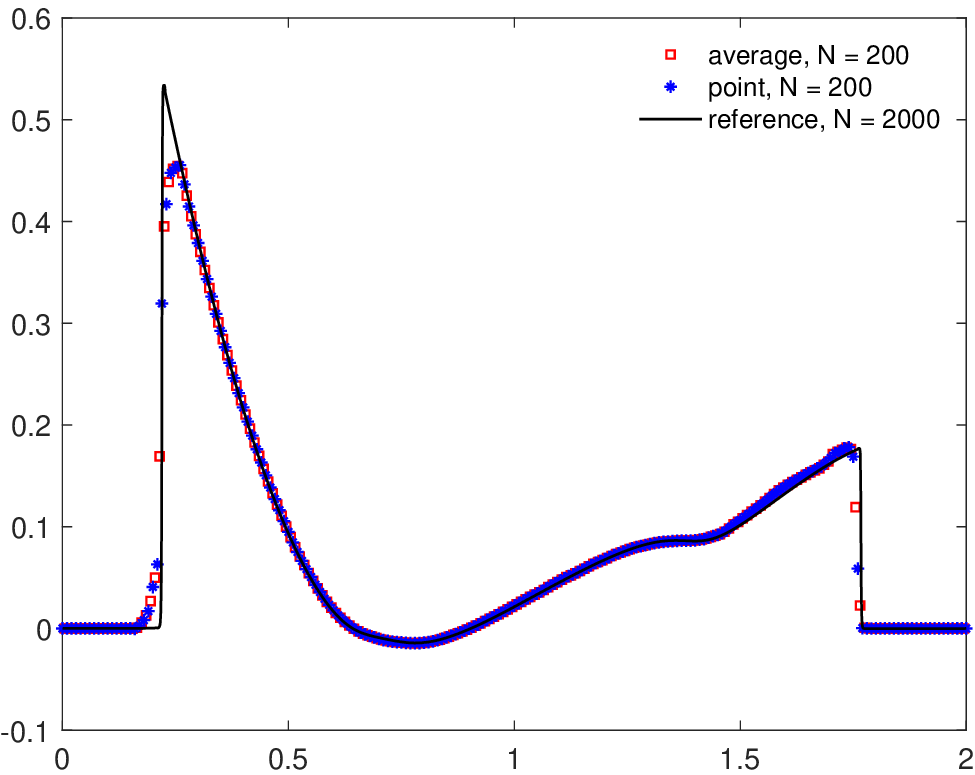}
  }\quad
  \subfigure[velocity ]{
  \centering
  \includegraphics[width=5.5cm,scale=1]{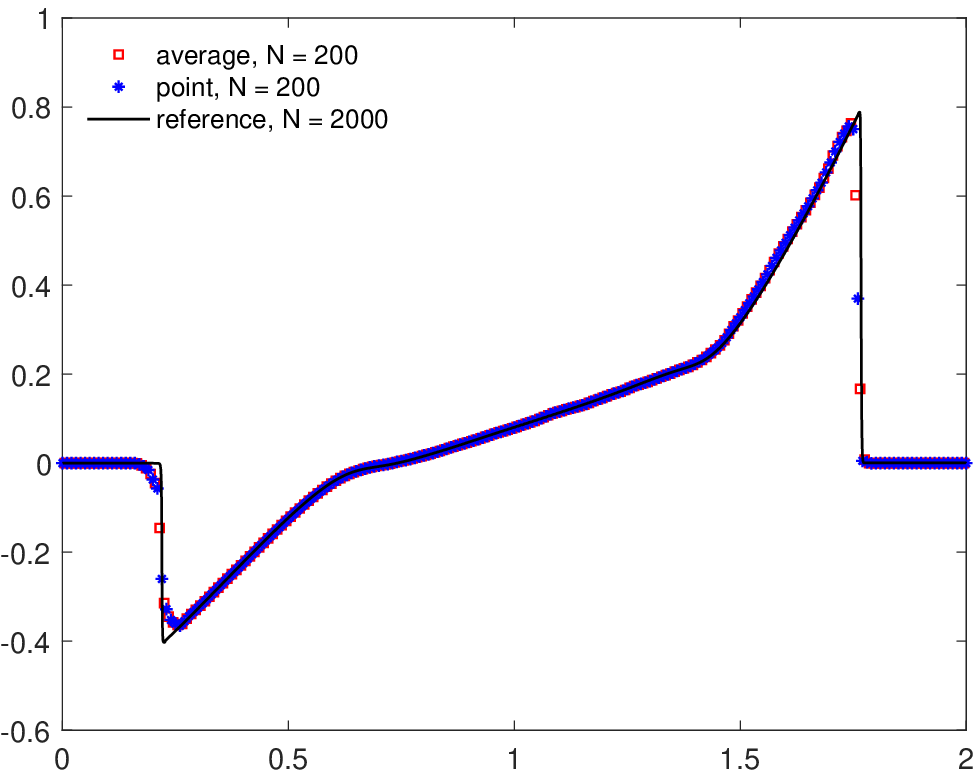}
  }
  \caption{Example \ref{euler:iseper_1d}: Pressure perturbation and velocity of the hydrostatic atmosphere for the large amplitude wave propagation $A=1$. }\label{fig:perbig_p}
\end{figure}

To further assess the performance of the proposed method for problems involving discontinuities, we consider the case with a large amplitude $A = 1$. The stopping time is kept the same as before. Fig. \ref{fig:perbig_p} shows the simulation results with 200 and 2000 uniform cells for comparison. It can be seen that our proposed method accurately captures the propagation of the perturbation while effectively controlling spurious oscillations near the discontinuity.

\begin{figure}[htb!]
  \centering
  \subfigure[density]{
  \centering
  \includegraphics[width=5.5cm,scale=1]{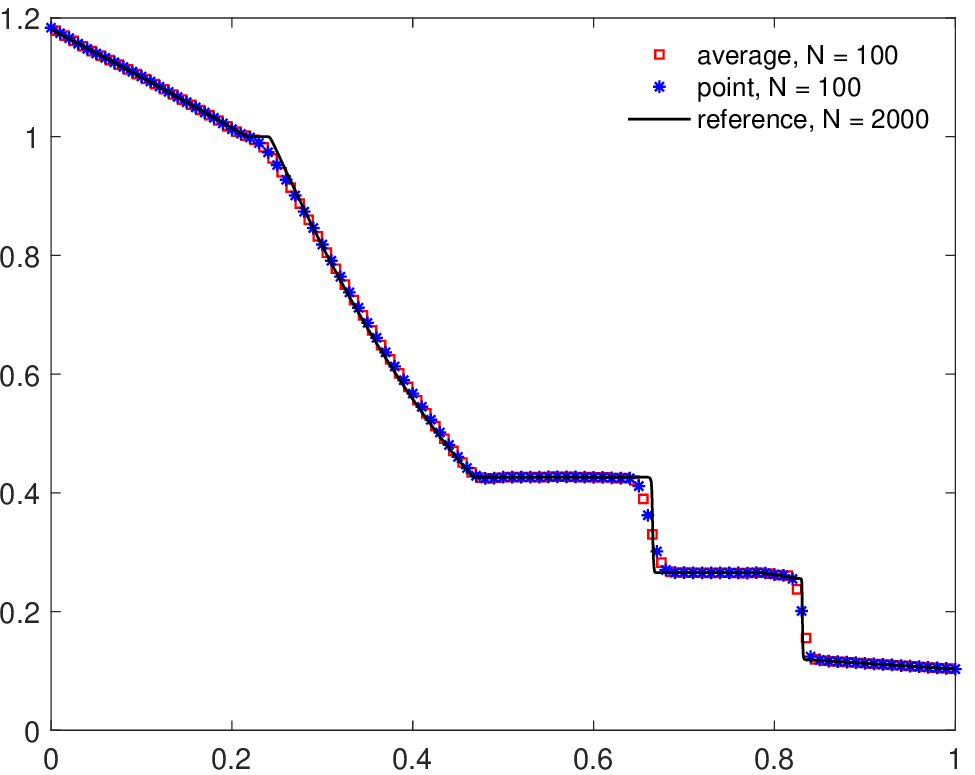}
  }\quad
  \subfigure[velocity]{
  \centering
  \includegraphics[width=5.5cm,scale=1]{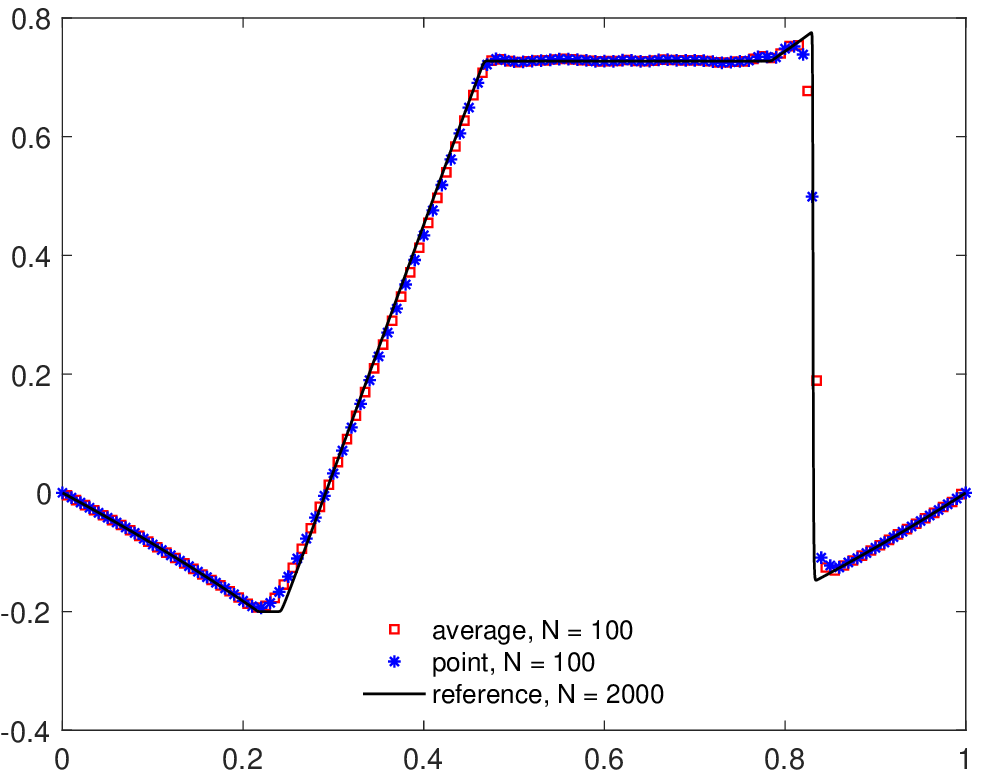}
  }

  \subfigure[energy]{
  \centering
  \includegraphics[width=5.5cm,scale=1]{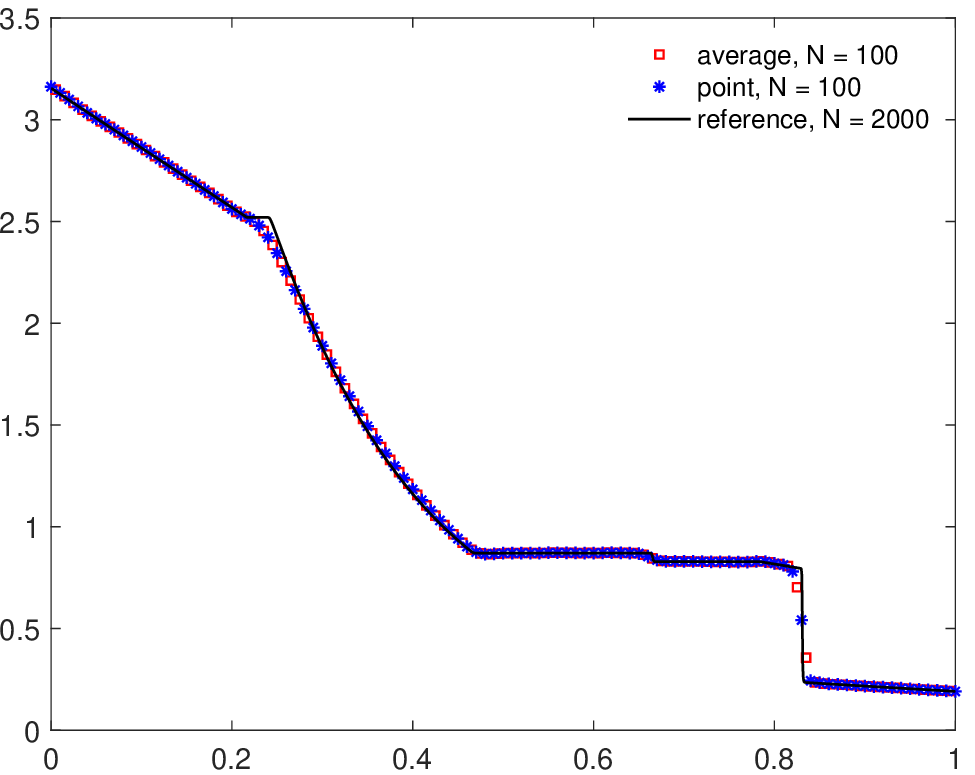}
  } \quad
  \subfigure[pressure]{
  \centering
  \includegraphics[width=5.5cm,scale=1]{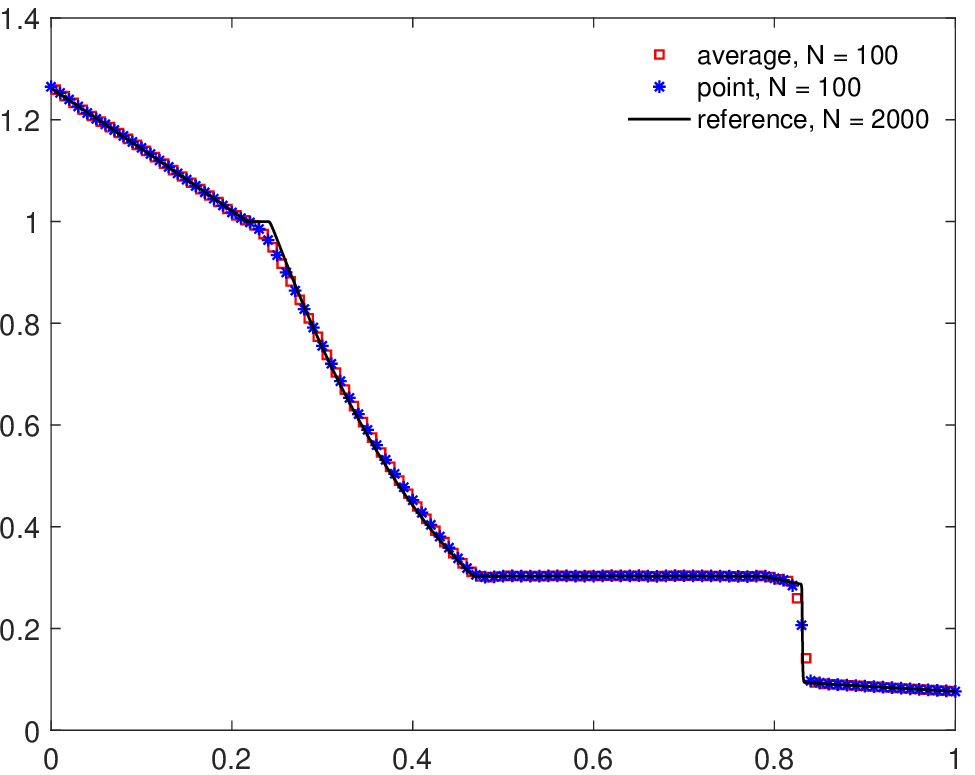}
  }
  \caption{Example \ref{euler:shock}: The numerical density, velocity, energy and pressure for the shock tube problem at time $t=0.2$.}\label{fig:shock}
\end{figure}

\begin{example}{\bf Shock tube problem}\label{euler:shock}
\end{example}
We consider the standard Sod test coupled with a linear gravitational field $\phi(x)=x$. The initial conditions are given by
\begin{equation}\label{euler:sod_ini}
  (\rho,u,p)(x,0) = \left\{\begin{array}{lll}
      (1,0,1),&  \text{if} \  x \leq 0.5, \\
      (0.125,0,0.1),&\text{otherwise},\\
      \end{array}\right.
\end{equation}
with $\gamma=1.4$ on the computational domain $[0,1]$. Solid wall boundary conditions are imposed, and the final time is set to $t=0.2$. Fig. \ref{fig:shock} shows the corresponding numerical results for the density, velocity, energy and pressure with a mesh of 100 cells. For comparison, the solutions obtained by the third-order well-balanced DG scheme \cite{li2018welldge} on 2000 uniform cells are treated as the reference solution. It is evident that our developed method achieves excellent agreement with the reference solution, and provides oscillation-free solutions.

\begin{figure}[htb!]
  \centering
  \subfigure[density]{
  \centering
  \includegraphics[width=4.9cm,scale=1]{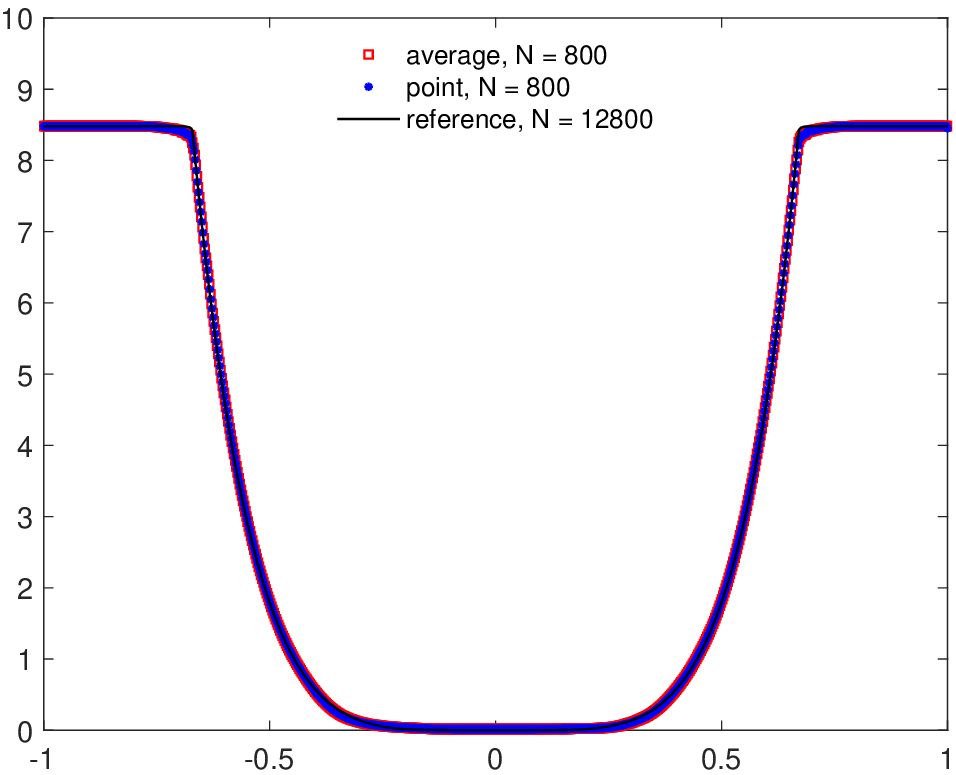}
  }
  \subfigure[momentum ]{
  \centering
  \includegraphics[width=4.9cm,scale=1]{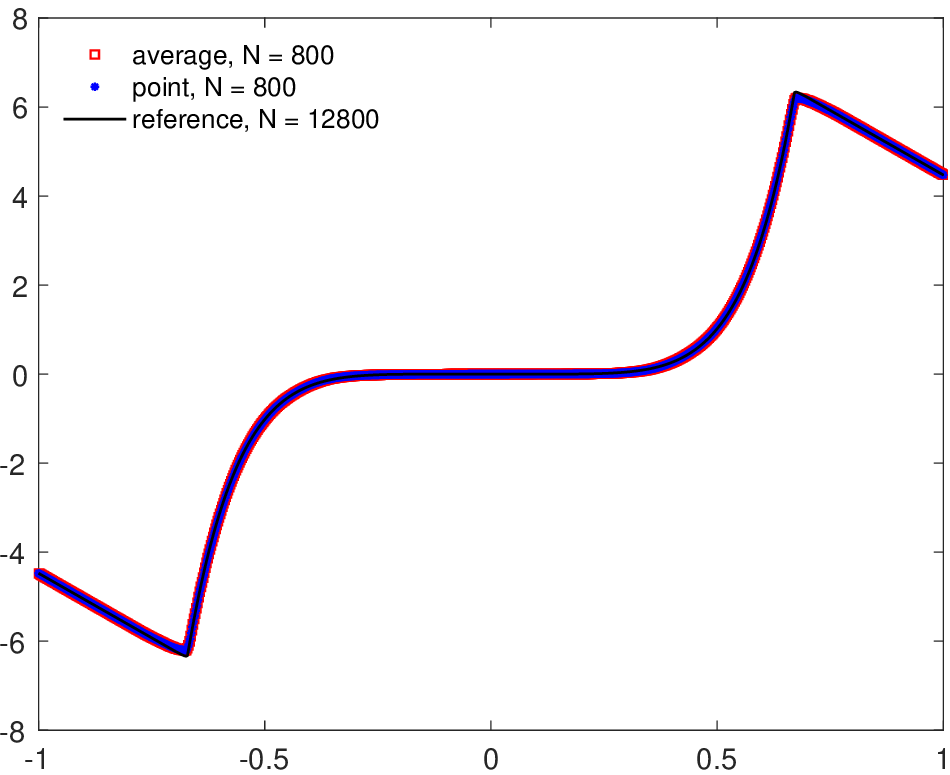}
  }
  \subfigure[energy]{
  \centering
  \includegraphics[width=4.9cm,scale=1]{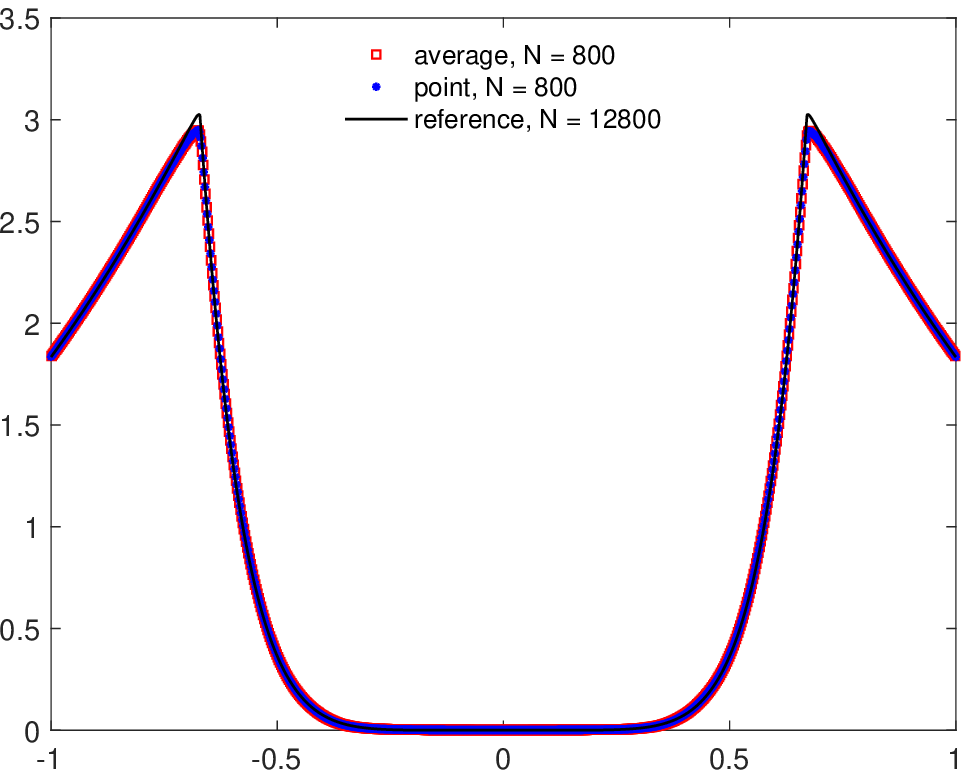}
  }
  \caption{Example \ref{euler:rare1d}: The numerical density, momentum and energy for the rarefaction test at time $t=0.6$.}\label{fig:rare}
\end{figure}

\begin{example}{\bf Rarefaction test with low density and low pressure}\label{euler:rare1d}
\end{example}
To test the positivity-preserving property of our constructed active flux method, we take an extreme rarefaction test on the computational domain $[-1,1]$. The initial state is the Riemann problem given by \cite{zhang2010positivity}, 
\begin{equation}\label{euler:rare_ini}
  \rho(x,0)=7, \quad p(x,0) = 0.2, \quad 
  u(x,0) = \left\{\begin{array}{lll}
    -1,&    x \leq 0, \\
    1,&x>0,\\
    \end{array}\right.
\end{equation}
coupled with a quadratic gravitational field $\phi(x) = x^2/2$ centered around $x=0$. Outflow boundary conditions are imposed at $x=-1$ and $x=1$. A MOOD paradigm is necessary since the presence of extremely low density and pressure. Without special positivity-preserving treatment, the scheme will blow up numerically in a few time steps. Fig. \ref{fig:rare} plots the numerical results at time $t=0.6$, calculated by our method with 800 and 12800 uniform cells for comparison. It is clear that the structure of the low density and low pressure are well captured with the MOOD paradigm, which indicates the effectiveness of our proposed structure-preserving active flux method.

\begin{figure}[htb!]
  \centering
  \subfigure[density]{
  \centering
  \includegraphics[width=4.9cm,scale=1]{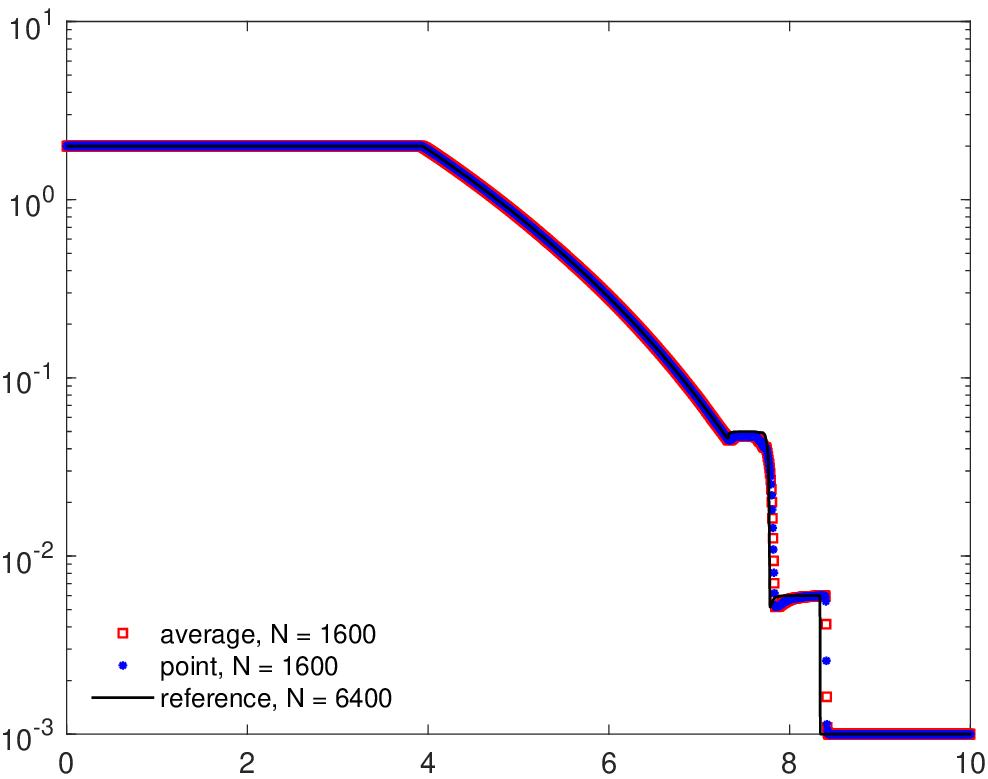}
  }
  \subfigure[momentum ]{
  \centering
  \includegraphics[width=4.9cm,scale=1]{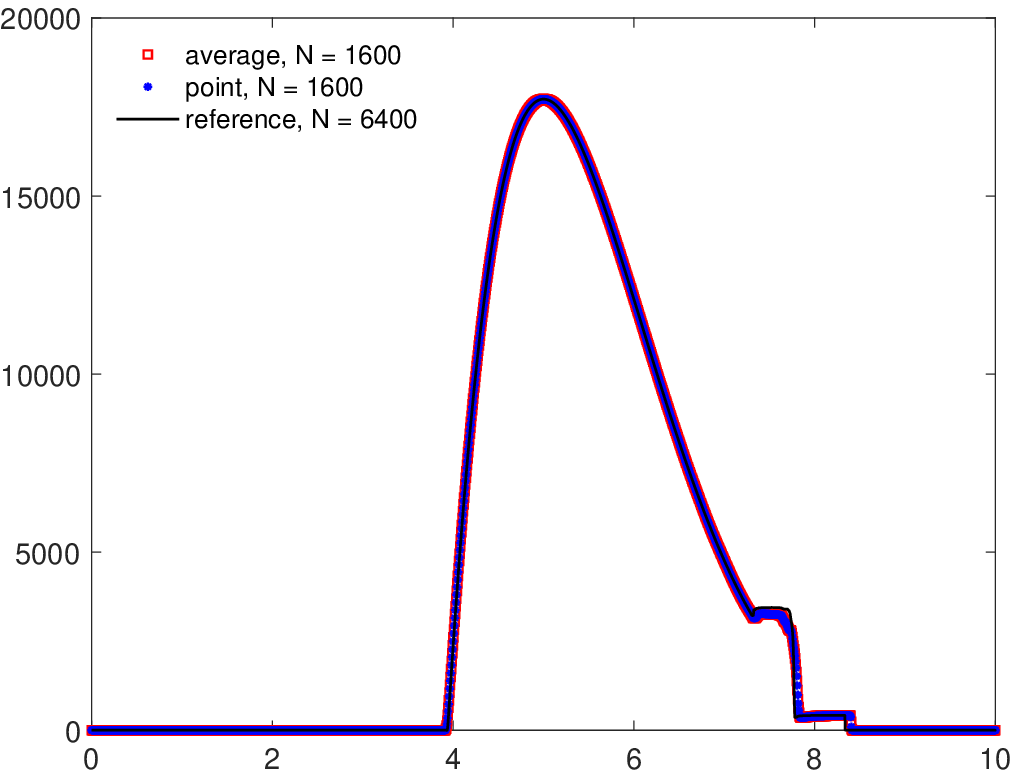}
  }
  \subfigure[pressure]{
  \centering
  \includegraphics[width=4.9cm,scale=1]{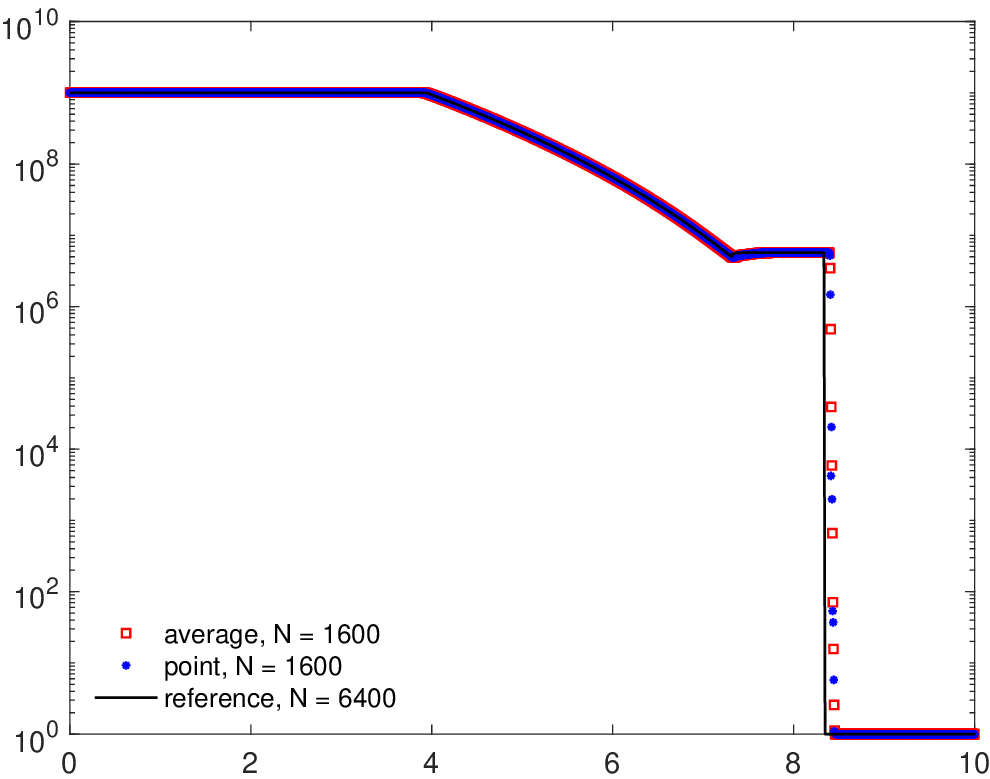}
  }
  \caption{Example \ref{euler:leblanc}: The log plot of density, momentum and the log plot of pressure for the leblanc problem at time $t=0.00004$.}\label{fig:leblanc}
\end{figure}

\begin{example}{\bf Leblanc problem in linear gravitational field}\label{euler:leblanc}
\end{example}
We consider another highly challenging test problem, which is an extension of one-dimensional Leblanc shock tube problem with linear gravitational potential $\phi(x) = x$. The initial conditions are given by  
\begin{equation}\label{euler:leblanc_ini}
  (\rho,u,p)(x,0) = \left\{\begin{array}{lll}
    (2,0,10^9),&   x \leq 5, \\
    (10^{-3},0,1),&x>5.\\
    \end{array}\right. 
\end{equation}
We can see that the density and pressure involve strong jumps initially. 
Solid wall boundary conditions are imposed at the end points $x=0$ and $x=10$.
A sufficiently fine mesh is required to fully resolve the wave structure in this example. To this end, we treat the solution obtained on a uniform mesh of 6400 cells as the reference solution. The numerical results at time $t=0.00004$ with 1600 uniform cells are also presented in Fig. \ref{fig:leblanc}. 
We can observe that the discontinuity is well captured without introducing spurious oscillations, while exhibiting excellent resolution of the wave structures.

\subsection{Two-dimensional tests}
\begin{example}{\bf Accuracy test}\label{euler:accuracy2D}
\end{example}
To verify the high-order accuracy of our structure-preserving active flux method for a smooth solution in two-dimensional Euler equations, we consider the following exact solutions on the computational domain $[0,2]\times[0,2]$,
\begin{equation}\label{euler:acc2d}
  \begin{aligned}
    & \rho(x,y,t) = 1 + 0.2\sin(\pi(x+y-t(u_0+v_0))), \\
    & u(x,y,t) = u_0, \\
    & v(x,y,t) = v_0, \\
    & p(x,y,t) = p_0 + t(u_0+v_0) - x - y + 0.2\cos(\pi(x+y-t(u_0+v_0)))/\pi,
  \end{aligned}
\end{equation}
where the parameters $u_0=v_0=1,p_0=4.5$ and $\gamma=1.4$. The gravity is given by a linear gravitational potential $\phi_x = \phi_y=1$. The boundary condition is specified by the exact solution. The final time is set as $t=0.1$.  The $L^1$ errors of the cell average and point values, along with their corresponding convergence rates are presented in Tables \ref{table:DGacc2D_ave} and \ref{table:DGacc2D_cell}, respectively. It is observed that our proposed structure-preserving active flux method attains third-order accuracy. 
\begin{table}[htb]
  \centering
  \caption{Example \ref{euler:accuracy2D}:  $L^1$ errors of the cell average and numerical orders of accuracy with initial condition (\ref{euler:acc2d}), $t = 0.1$. }\label{table:DGacc2D_ave}
  \begin{tabular}{c c c c c c c c c c}
  \toprule
    \multirow{2}{*} {$Nx\times Ny$} & \multicolumn{2}{c}{$\rho$}&\multicolumn{2}{c}{$ \rho u $}&\multicolumn{2}{c}{$ \rho v$}& \multicolumn{2}{c}{$E$}\\
   \cmidrule(lr){2-3} \cmidrule(lr){4-5} \cmidrule(lr){6-7}\cmidrule(lr){8-9}
    ~ &$L^1$ error &order&$L^1$ error &order&$L^1$ error &order&$L^1$ error &order\\
   \midrule
   10$\times$10 &  2.68E-04 & -- & 2.95E-04 & -- & 2.95E-04 & -- & 9.32E-04 & -- \\
   20$\times$20  & 2.73E-05 & 3.30 & 3.37E-05 & 3.13 & 3.37E-05 & 3.13 & 1.34E-04 & 2.79 \\
   40$\times$40 & 3.47E-06 & 2.98 & 4.46E-06 & 2.92 & 4.46E-06 & 2.92 & 1.80E-05 & 2.90 \\
   80$\times$80 & 4.52E-07 & 2.94 & 5.81E-07 & 2.94 & 5.81E-07 & 2.94 & 2.33E-06 & 2.95 \\
   160$\times$160 & 5.80E-08 & 2.96 & 7.45E-08 & 2.96 & 7.45E-08 & 2.96 & 2.95E-07 & 2.98\\
  \bottomrule
  \end{tabular}
\end{table}

\begin{table}[htb]
  \centering
  \caption{Example \ref{euler:accuracy2D}:  $L^1$ errors of the point values and numerical orders of accuracy with initial condition (\ref{euler:acc2d}), $t = 0.1$. }\label{table:DGacc2D_cell}
  \begin{tabular}{c c c c c c c c c c}
  \toprule
    \multirow{2}{*} {$Nx\times Ny$} & \multicolumn{2}{c}{$\rho$}&\multicolumn{2}{c}{$ \rho u $}&\multicolumn{2}{c}{$ \rho v$}& \multicolumn{2}{c}{$E$}\\
   \cmidrule(lr){2-3} \cmidrule(lr){4-5} \cmidrule(lr){6-7}\cmidrule(lr){8-9}
    ~ &$L^1$ error &order&$L^1$ error &order&$L^1$ error &order&$L^1$ error &order\\
   \midrule
   10$\times$10 & 1.48E-03 & -- & 1.42E-03 & --& 1.42E-03 & -- & 2.38E-03 & -- \\
   20$\times$20  & 1.32E-04 & 3.48 & 1.32E-04 & 3.43 & 1.32E-04 & 3.43 & 2.67E-04 & 3.16 \\
   40$\times$40  & 1.49E-05 & 3.15 & 1.52E-05 & 3.11 & 1.52E-05 & 3.11 & 3.26E-05 & 3.03 \\
   80$\times$80  & 1.82E-06 & 3.03 & 1.89E-06 & 3.01 & 1.89E-06 & 3.01 & 4.07E-06 & 3.00 \\
   160$\times$160 & 2.27E-07 & 3.00 & 2.37E-07 & 2.99 & 2.37E-07 & 2.99 & 5.09E-07 & 3.00\\
  \bottomrule
  \end{tabular}
\end{table}

\begin{example}{\bf Test for two-dimensional isothermal equilibrium}\label{euler:Poly2_2d}
\end{example}
To illustrate the capability of our developed method for the preservation and perturbation of the isothermal hydrostatic equilibrium state in two dimensions, we test the same example as in \cite{li2016highfve}, where the isothermal equilibrium state is given by
\begin{equation}\label{euler:isepoly_2d}
  \begin{aligned}
  &\rho(x,y) = \rho_0\exp\left(-\dfrac{\rho_0 g}{p_0}(x+y)\right), \\
  & u(x,y)=v(x,y) = 0,\\
  & p(x,y) = p_0\exp\left(-\dfrac{\rho_0 g}{p_0}(x+y)\right),
  \end{aligned}
\end{equation}
where the parameters $\rho_0=1.21,p_0=1,g=1$. 
The linear gravitation potential $ \phi(x,y) = x+y$ is considered on the  computational domain $[0,1]\times[0,1]$. The ratio of specific heats is $\gamma=1.4$. We simulate the solutions up to $t = 0.5$ with different mesh sizes and present the $L^1$ and $L^{\infty}$ errors of the cell average and point values at double precision in Table \ref{table:isewb_1002D}.  It can be seen clearly that the isothermal equilibrium state is exactly preserved even on coarse meshes in two-dimensional case.

\begin{table}[htb]
  \centering
  \caption{Example \ref{euler:Poly2_2d}: $L^1$ and $L^{\infty}$ errors for two-dimensional isothermal hydrostatic equilibrium flow (\ref{euler:isepoly_2d}).}\label{table:isewb_1002D}
  \begin{tabular}{c c c c c c c }
   \toprule
   ~&{$Nx\times Ny$}&{error}&{$\rho$}&{$\rho u$}&{$\rho v$}&{$E$}\\
    \midrule
    \multirow{4}{*} {average} &\multirow{2}{*} {50$\times$50}  &$L^1$   & 1.98E-16 & 2.58E-16 & 2.38E-16 & 7.49E-16 \\
    ~& &$L^{\infty}$  &3.77E-15 & 6.23E-15 & 5.92E-15 & 1.42E-14 \\
    \cmidrule(lr){2-7}
    &\multirow{2}{*} {100$\times$100}  &$L^1$   &2.61E-16 & 1.19E-16 & 1.31E-16 & 3.98E-16 \\
    ~& &$L^{\infty}$  & 7.99E-15 & 2.21E-15 & 2.14E-15 & 6.22E-15\\
    \midrule
    \multirow{4}{*} {point } &\multirow{2}{*} {50$\times$50}  &$L^1$   &  2.17E-16 & 2.86E-16 & 2.67E-16 & 8.09E-16 \\
    ~& &$L^{\infty}$  &4.22E-15 & 6.94E-15 & 6.41E-15 & 1.42E-14 \\
    \cmidrule(lr){2-7}
    &\multirow{2}{*} {100$\times$100}  &$L^1$   &2.61E-16 & 1.27E-16 & 1.39E-16 & 4.13E-16 \\
    ~& &$L^{\infty}$  & 5.22E-15 & 2.48E-15 & 2.36E-15 & 6.22E-15\\
    \bottomrule
  \end{tabular}
\end{table}

\begin{figure}[htb!]
  \centering
  \subfigure[pressure perturbation]{
  \centering
  \includegraphics[width=5.5cm,scale=1]{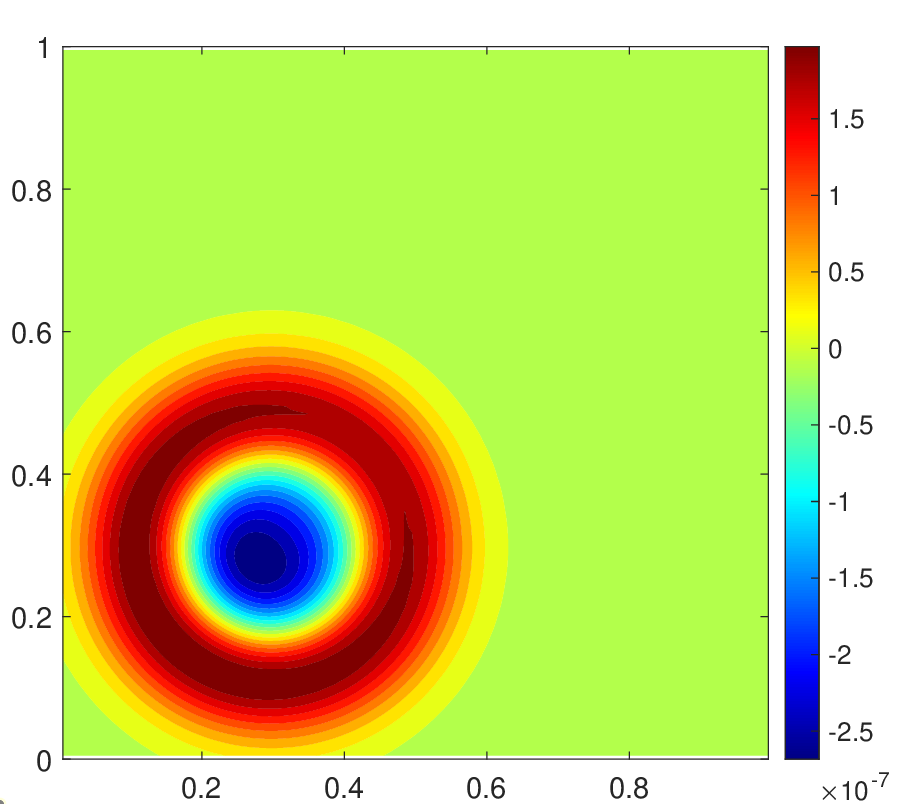}
  }\quad
  \subfigure[velocity]{
  \centering
  \includegraphics[width=5.5cm,scale=1]{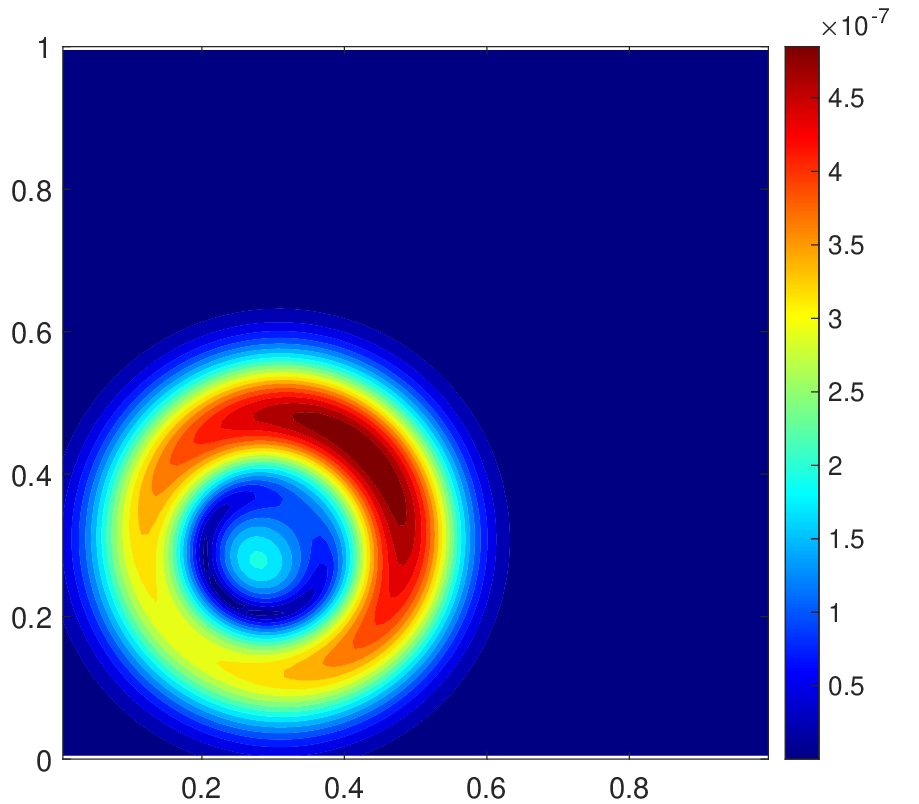}
  }
  \subfigure[pressure perturbation]{
  \centering
  \includegraphics[width=5.5cm,scale=1]{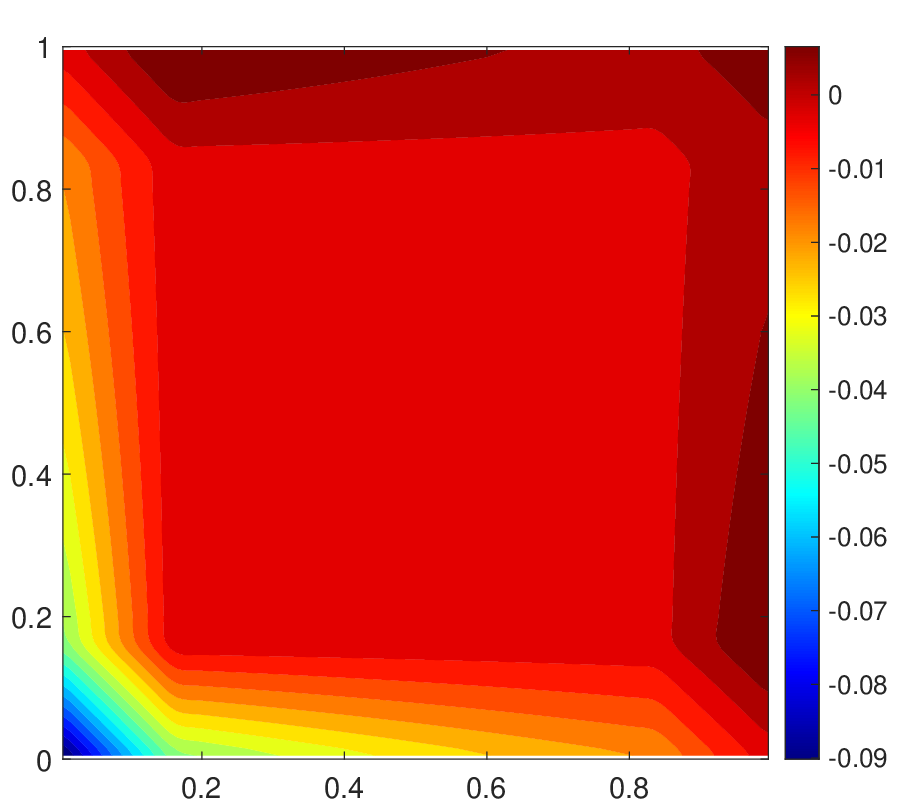}
  }\quad
  \subfigure[velocity]{
  \centering
  \includegraphics[width=5.5cm,scale=1]{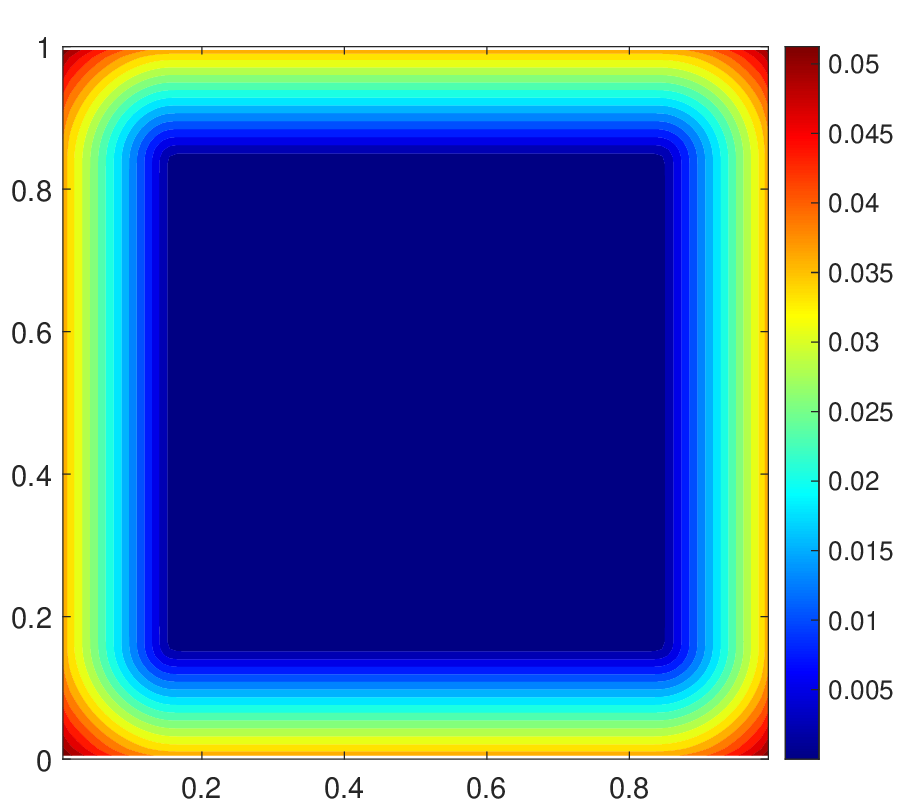}
  }
  \caption{Example \ref{euler:Poly2_2d}: The cell average contours of the pressure perturbation and velocity $\sqrt{u^2+v^2}$ for the small amplitude $A=10^{-6}$ in (\ref{euler:ise_ppres_2d}), with 20 uniformly spaced contour lines at time $t=0.15$. Top: results by our well-balanced active flux method. Bottom: results by the non-well-balanced method. }\label{fig:hypoly2_wbnwb2d}
\end{figure}

Subsequently, a small perturbation is added to the initial pressure profile
 \begin{equation}\label{euler:ise_ppres_2d}
  \begin{aligned}
  p(x,y,0) = p(x,y) + A\exp\left(-\frac{100\rho_0}{p_0}((x-0.3)^2+(y-0.3)^2)\right),
  \end{aligned}
\end{equation}
with the density and velocities remaining unchanged. We consider the small amplitude $A=10^{-6}$ and impose transmission boundary conditions. We simulate the solution up to $t = 0.15$ on a uniform grid of 100$\times$100 cells, and display the cell average contours of the pressure perturbation and velocity in Fig. \ref{fig:hypoly2_wbnwb2d}. For comparison, the corresponding results obtained by the non-well-balanced method are also illustrated. The results indicate that our proposed well-balanced scheme outperforms the non-well-balanced counterpart. It can achieve high resolution and is capable of capturing the small features of the disturbance accurately.

\begin{example}{\bf  Test for two-dimensional polytropic equilibrium}\label{euler:Poly1_2d}
\end{example}
We use this example to assess the performance of our developed scheme near the two-dimensional polytropic equilibrium. An astrophysical problem studied in\cite{kappeli2014well} is considered on the computational domain $[-0.5,0.5]\times[-0.5,0.5]$. The analytical hydrostatic equilibrium and the gravitational potential are given by
\begin{equation}\label{euler:hypoly_2d}
  \begin{aligned}
  &\rho(r) = \rho_c\dfrac{\sin(\alpha r)}{\alpha r}, \ u(r) =v(r) = 0, \ p(r) = \rho(r)^2, \\
  &\phi(r) = -2\rho_c\dfrac{\sin(\alpha r)}{\alpha r},
  \end{aligned}
\end{equation}
where $r = \sqrt{x^2+y^2}$ is the radial variable, the parameters are taken as $\alpha = \sqrt{2\pi}$, $\rho_c=1$ and $\gamma=2$. To demonstrate that the steady state is indeed preserved up to machine error, we calculate the $L^1$ and $L^{\infty}$ errors at double precision up to $t = 2$ with different mesh sizes. The corresponding results both for the cell average and point values are listed in Table \ref{table:polywb_1002D}, from which the desired well-balanced property can be observed.
\begin{table}[htb]
  \centering
  \caption{Example \ref{euler:Poly1_2d}: $L^1$ and $L^{\infty}$ errors for two-dimensional polytropic hydrostatic equilibrium flow (\ref{euler:hypoly_2d}). }\label{table:polywb_1002D}
  \begin{tabular}{c c c c c c c }
   \toprule
    &{$Nx\times Ny$}&{error}&{$\rho$}&{$\rho u$}&{$\rho v$}&{$E$}\\
    \midrule
    \multirow{4}{*} {average}&\multirow{2}{*} {50$\times$50} &$L^1$        & 9.89E-17 & 1.59E-16 & 1.60E-16 & 6.86E-17 \\
    &~&$L^{\infty}$ & 1.44E-15 & 1.05E-15 & 1.36E-15 & 5.00E-16 \\
    \cmidrule(lr){2-7}
    &\multirow{2}{*} {100$\times$100} &$L^1$        &1.53E-16 & 2.27E-16 & 2.27E-16 & 7.21E-17 \\
    &~&$L^{\infty}$ &3.00E-15 & 1.45E-15 & 1.49E-15 & 6.66E-16\\
    \midrule
    \multirow{4}{*} {point}&\multirow{2}{*} {50$\times$50} &$L^1$        & 1.38E-16 & 1.79E-16 & 1.80E-16 & 8.23E-17 \\
    &~&$L^{\infty}$ &  8.88E-16 & 1.12E-15 & 1.17E-15 & 4.44E-16 \\
    \cmidrule(lr){2-7}
    &\multirow{2}{*} {100$\times$100} &$L^1$        &1.79E-16 & 2.45E-16 & 2.45E-16 & 8.18E-17 \\
    &~&$L^{\infty}$ &2.11E-15 & 1.68E-15 & 1.75E-15 & 6.66E-16\\
    \bottomrule
  \end{tabular}
\end{table}

\begin{figure}[htb!]
  \centering
  \subfigure[pressure perturbation]{
  \centering
  \includegraphics[width=5.5cm,scale=1]{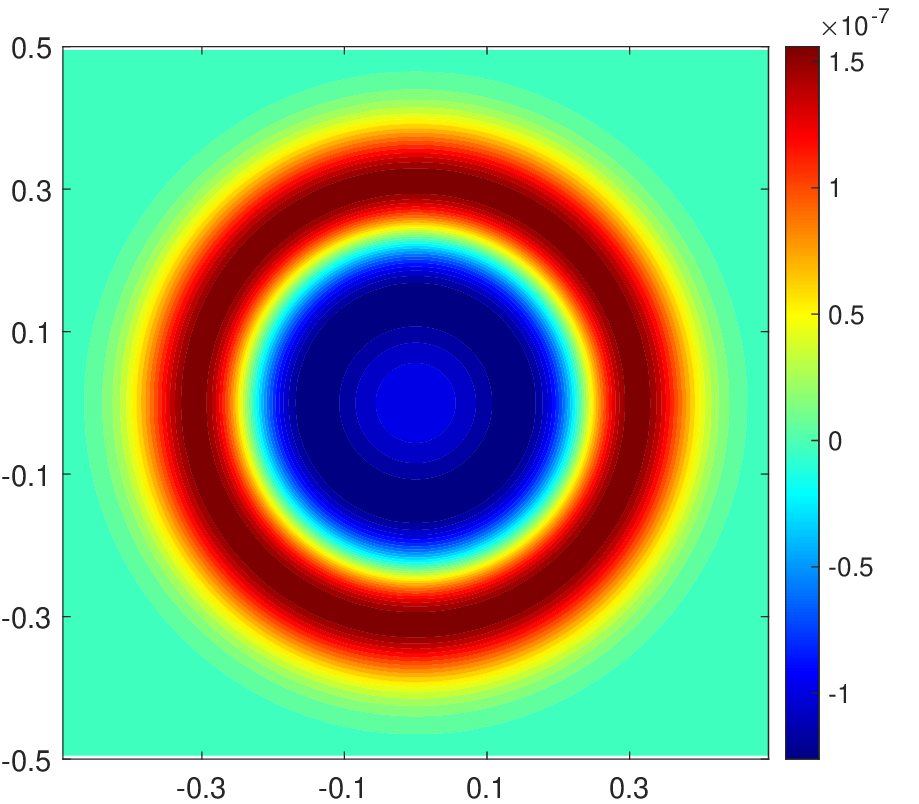}
  }\quad
  \subfigure[velocity]{
  \centering
  \includegraphics[width=5.5cm,scale=1]{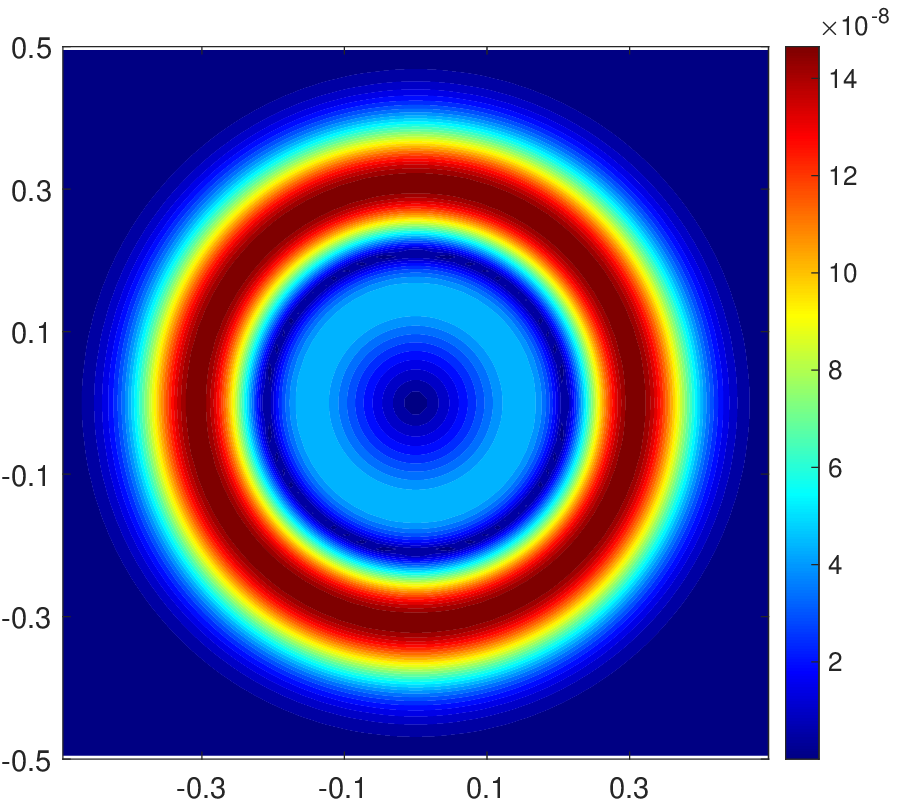}
  }
  \subfigure[pressure perturbation]{
  \centering
  \includegraphics[width=5.5cm,scale=1]{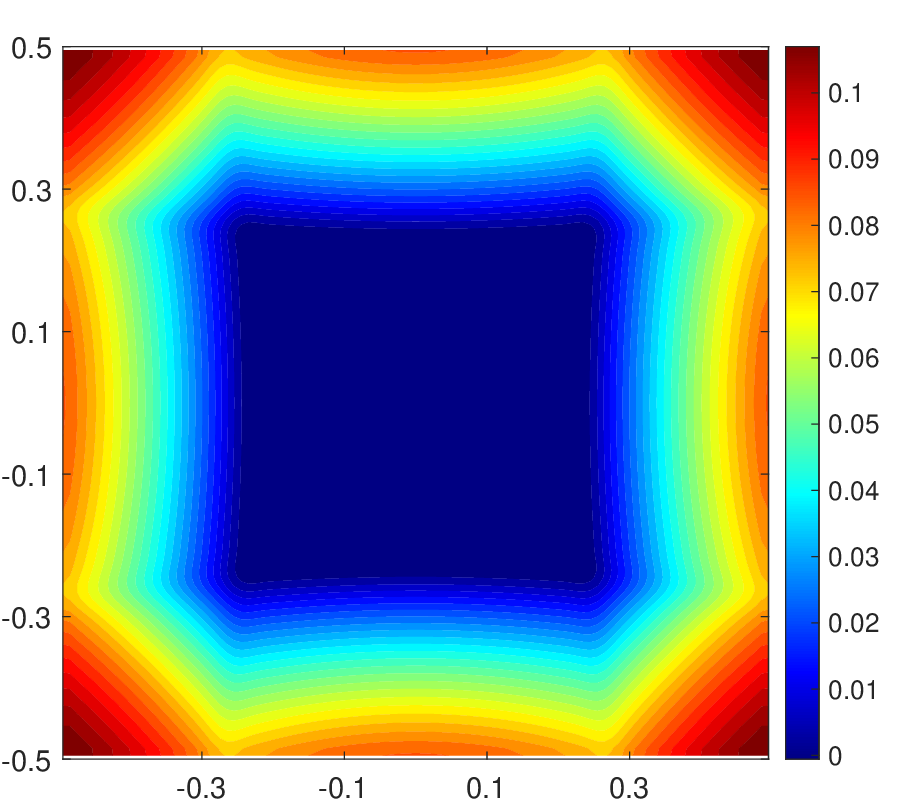}
  }\quad
  \subfigure[velocity]{
  \centering
  \includegraphics[width=5.5cm,scale=1]{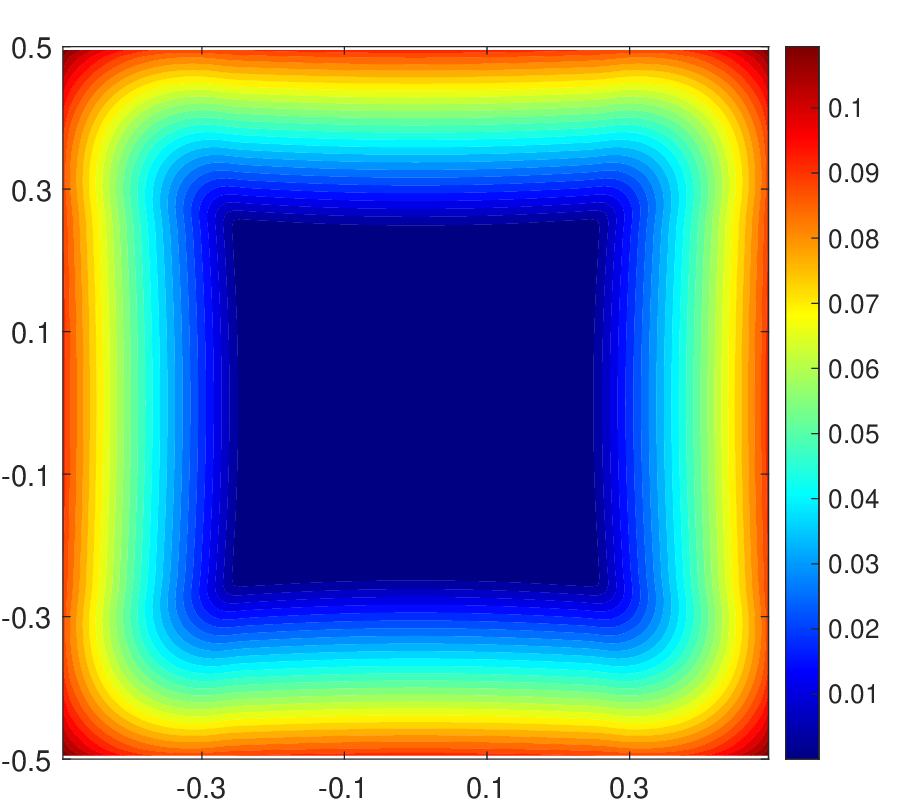}
  }
  \caption{Example \ref{euler:Poly1_2d}: The cell average contours of the pressure perturbation and velocity $\sqrt{u^2+v^2}$ for the small amplitude $A=10^{-6}$ in (\ref{euler:hyp_ppres_2d}), with 30 uniformly spaced contour lines at time $t=0.2$. Top: results by our well-balanced active flux method. Bottom: results by the non-well-balanced method. }\label{fig:hypoly1_wbnwb2d}
\end{figure}

We keep the density and velocity unchanged, and then add a small Gaussian hump perturbation to the initial pressure profile
\begin{equation}\label{euler:hyp_ppres_2d}
  \begin{aligned}
  p(x,y,0) = \rho(r)^2 + A\exp(-100r^2),
  \end{aligned}
\end{equation}
with a small amplitude $A=10^{-6}$. We impose transmission boundary conditions and simulate the solution to the final time $t = 0.2$ before the excited waves reach the boundary. The contours of the pressure perturbation and velocity, computed using our proposed well-balanced active flux method and a non-well-balanced method, are shown in Fig. \ref{fig:hypoly1_wbnwb2d} for comparison. The non-well-balanced method yields an inaccurate solution, while our proposed method performs robustly and maintains axial symmetry. This underscores the importance of the well-balanced scheme for accurately capturing small perturbations near the steady state.

\begin{figure}[htb!]
  \centering
  \subfigure[density]{
  \centering
  \includegraphics[width=4.9cm,scale=1]{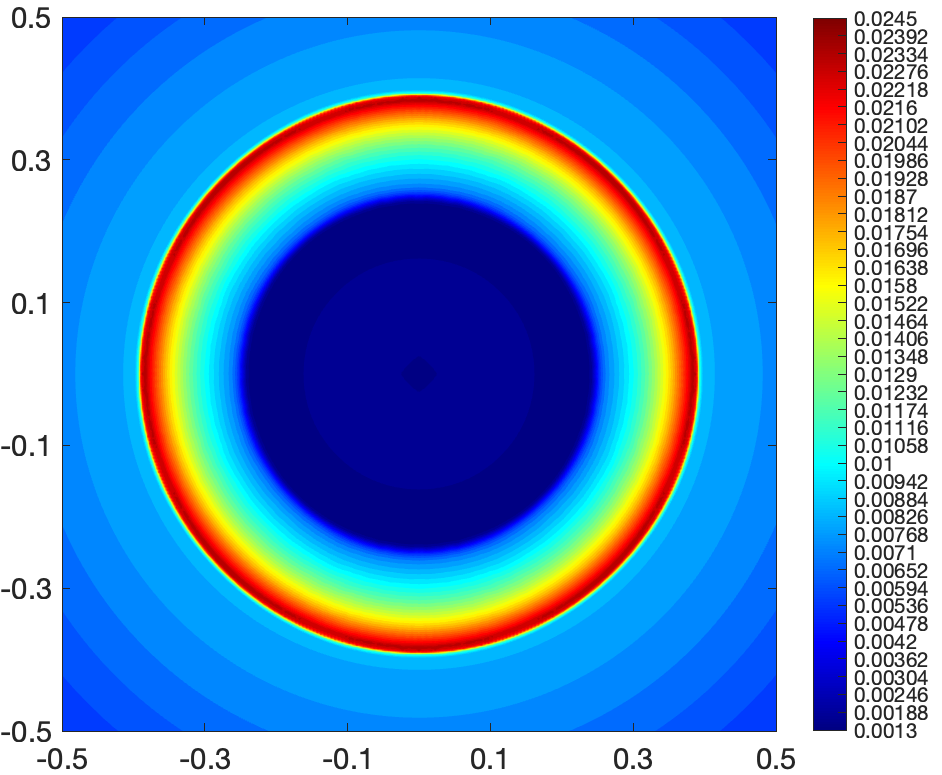}
  }
  \subfigure[pressure logarithm]{
  \centering
  \includegraphics[width=4.9cm,scale=1]{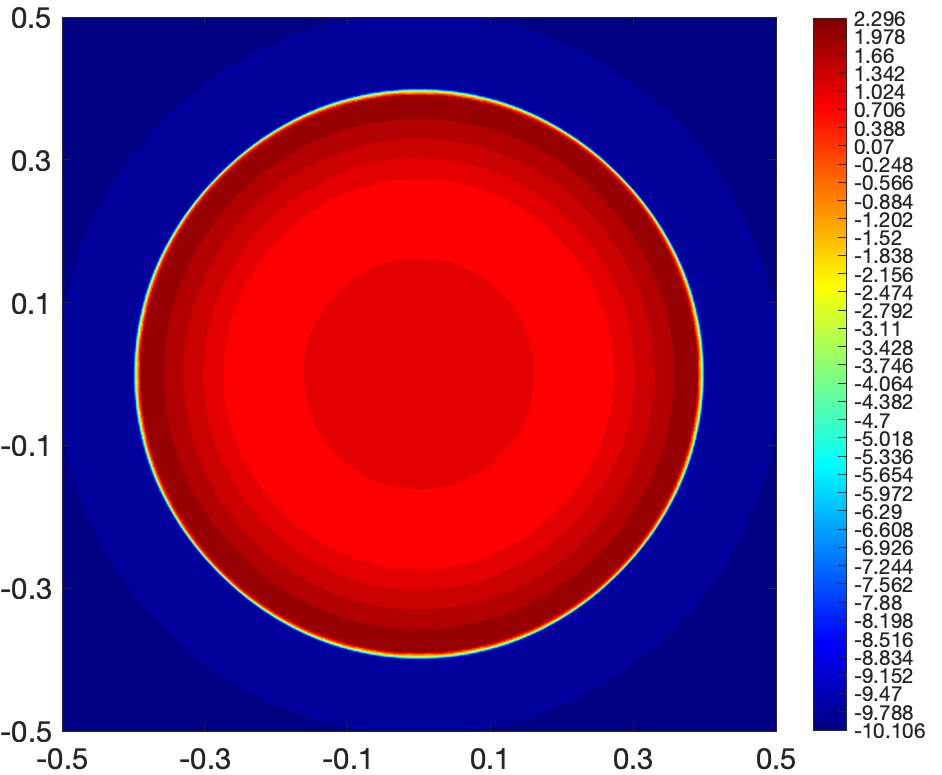}
  }
  \subfigure[pressure along $y=x$]{
  \centering
  \includegraphics[width=4.9cm,scale=1]{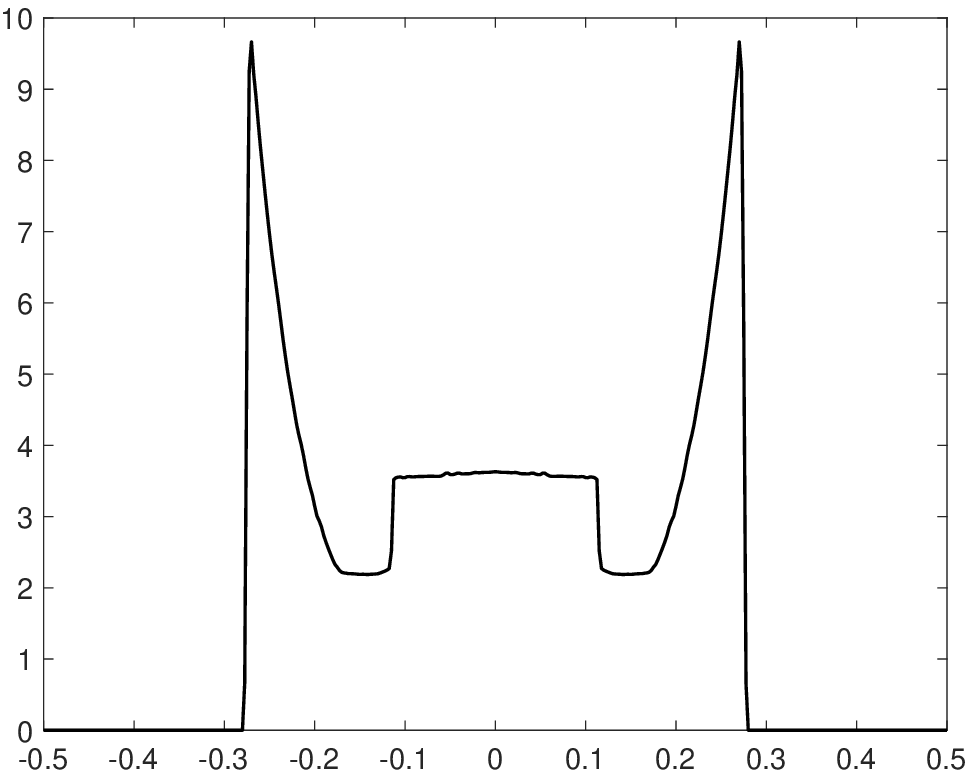}
  }
  \caption{Example \ref{euler:blast2d}: The point values contours of the density, the pressure logarithm, and the plot of pressure along $y=x$ at time $t=0.0005$. }\label{fig:blast}
\end{figure}

\begin{example}{\bf Two-dimensional blast problem}\label{euler:blast2d}
\end{example}
To further illustrate the positivity-preserving property of our scheme, we test a two-dimensional blast problem, which is introduced in \cite{wu2021uniformly}. The gravitational field, density, and velocity are chosen to be the same as those in \eqref{euler:hypoly_2d}, we only introduce a huge jump in the initial pressure, that is 
\begin{equation}\label{blast_per}
  p(x,y,0) = \rho(r)^2 + \left\{\begin{array}{lll}
    100,& r < 0.1, \\
    0,  & r\geq 0.1.\\
    \end{array}\right.
\end{equation}
The same set of parameters as in Example \ref{euler:Poly1_2d} is taken, except that $\rho_c$ is changed to 0.01, yielding the appearance of low density and low pressure in the solution. Transmissive boundary conditions are considered here. The computation is carried out on a 400$\times$400 mesh up to the final time $t=0.0005$, and we show the contours of the numerical density and pressure logarithm in Fig. \ref{fig:blast}. We also display the pressure profiles along the line $y=x$. It shows that the solution develops strong discontinuities, with a shock located at $\sqrt{x^2+y^2}/\sqrt{2}\approx0.28$. With the MOOD paradigm, our method is successful in maintaining the positivity-preserving property and suppressing numerical oscillations at the same time.

\begin{example}{\bf Inertia-gravity wave}\label{euler:ig}
\end{example}
The last example we take is given by \cite{ghosh2016well}, which is a benchmark problem arising from atmospheric flows. We consider a two-dimensional channel defined on the domain $[0,300000]\times[0,10000] \ \mathrm{m}^2$, with periodic boundary conditions on the left and right boundaries, and solid wall boundary conditions on the bottom and top boundaries. The flow evolves under a linear gravitational field with $\phi_x=0$ and $\phi_y=g=9.8\ \mathrm{m/s^2}$. The initial state consists of a uniformly stratified atmosphere with a constant velocity field $\boldsymbol{u}=(20\ \mathrm{m/s},0\ \mathrm{m/s})$. The potential temperature and Exner pressure are specified by
\begin{equation}
  \Theta = T_0\exp\left(\frac{\mathcal{N}^2}{g}y\right),\quad \Pi = 1+\frac{(\gamma-1)g^2}{\gamma RT_0 \mathcal{N}^2}\left[\exp\left(-\frac{\mathcal{N}^2}{g}y\right)-1\right].
\end{equation}
where the Brunt--V{\"a}is{\"a}l{\"a} frequency $\mathcal{N}=0.01/\mathrm{s}$, the gas constant $R=287.058\ \mathrm{J/kg\ K}$, and the reference temperature $T_0=300\ \mathrm{K}$ at $y=0\ \mathrm{m}$. A small perturbation is introduced into the initial potential temperature,
\begin{equation*}
  \Delta \Theta(x,y,0) = \theta_c \sin\left(\frac{\pi y}{h_c}\right)\left[1+{(x-x_c)^2}/{a_c^2}\right]^{-1},
\end{equation*}
with $\theta_c=0.01\ \mathrm{K}$, $h_c = 10000\ \mathrm{m}$, $x_c=10000\ \mathrm{m}$, $a_c=5000\ \mathrm{m}$. The pressure and density are obtained from $\Theta$ and $\Pi$ through
\begin{equation*}
  p = p_0\Pi^{\frac{\gamma}{\gamma-1}}, \quad \rho = \frac{p_0}{R\Theta}\Pi^{\frac{1}{\gamma-1}},
\end{equation*}
where the reference pressure is taken as $p_0=10^5\ \mathrm{N/m^2}$ at $y=0\ \mathrm{m}$.

\begin{figure}[htb!]
  \centering
  \subfigure[$t=0$]{
  \centering
  \includegraphics[width=7.5cm,scale=1]{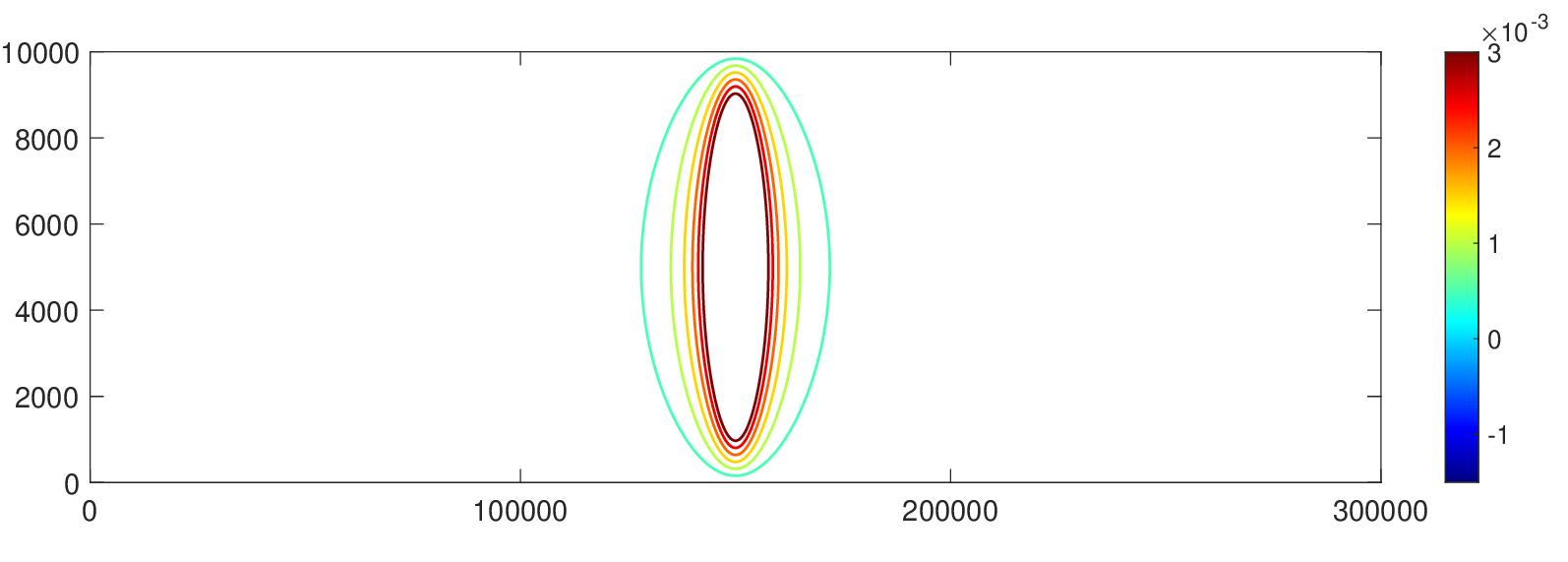}
  }
  \subfigure[$t=1000$]{
  \centering
  \includegraphics[width=7.5cm,scale=1]{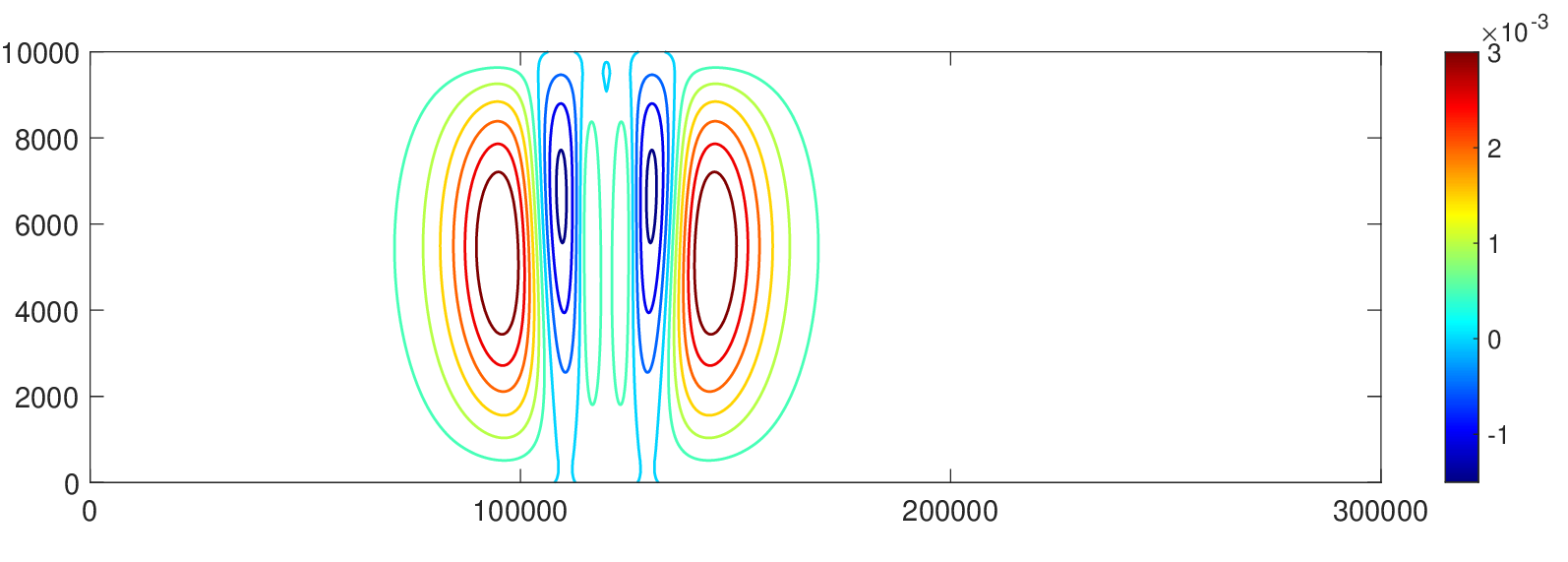}
  }
  \subfigure[$t=2000$]{
  \centering
  \includegraphics[width=7.5cm,scale=1]{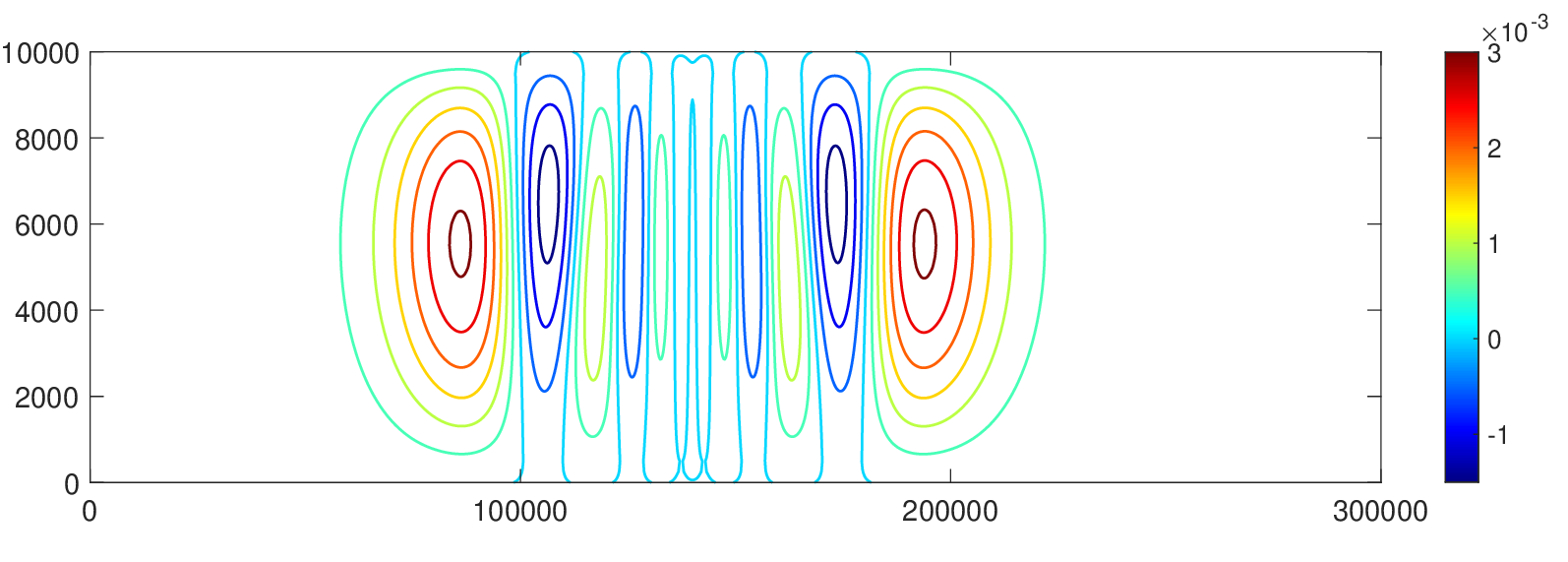}
  }
  \subfigure[$t=3000$]{
  \centering
  \includegraphics[width=7.5cm,scale=1]{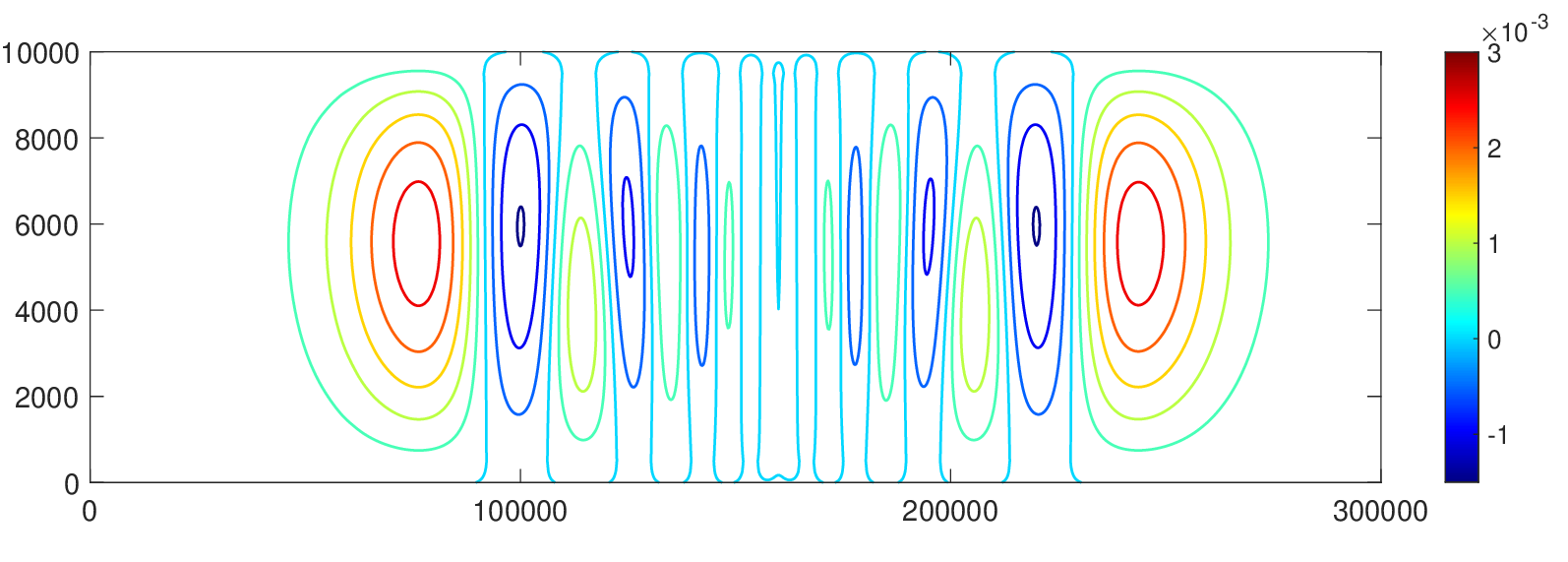}
  }
  \caption{Example \ref{euler:ig}: The point values contours of the potential temperature perturbation $\Delta \Theta$ at different times $t=0, 1000, 2000, 3000\ \mathrm{s}$. 10 uniformly spaced contour lines from $-0.0015$ to 0.003. }\label{fig:inertia}
\end{figure}

\begin{figure}[htb!]
  \centering
  \includegraphics[width=7cm,scale=1]{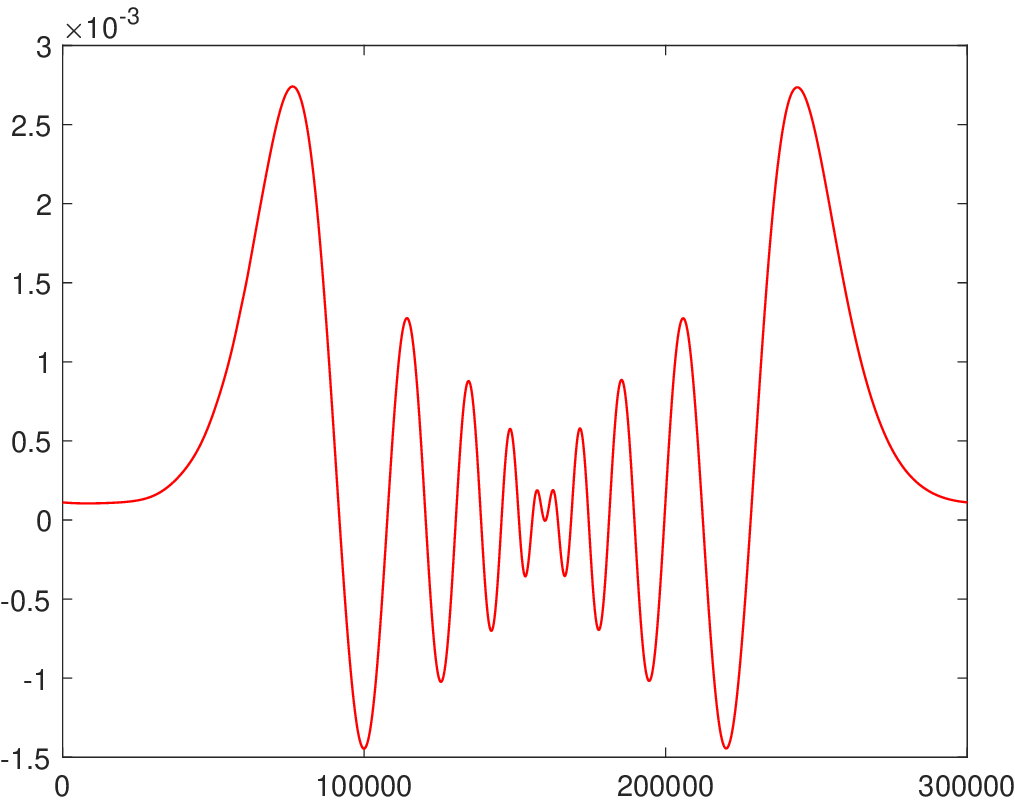}
  \caption{Example \ref{euler:ig}: The potential temperature perturbation $\Delta \Theta$ at time $t=3000\ \mathrm{s}$ along the line $y=5000\ \mathrm{m}$. }\label{fig:inertiam}
\end{figure}

We run the simulation until $t=3000\ \mathrm{s}$ with 1200$\times$40 mesh cells, seeing the contours of the potential temperature perturbation $\Delta \Theta$ at different times $t=0, 1000, 2000, 3000\ \mathrm{s}$ in Fig. \ref{fig:inertia}, and the profiles of $\Delta \Theta$ at $t=3000\ \mathrm{s}$ along the line $y=5000\ \mathrm{m}$ in Fig. \ref{fig:inertiam}. We can observe a well-resolved evolution of the potential temperature perturbation, which demonstrates excellent agreement with the reference structures in \cite{ghosh2016well,wu2021uniformly}.

\section{Conclusion}\label{se:co}
In this work, the high-order structure-preserving active flux method was constructed for the compressible Euler equations under gravitational fields. Accounting for the well-balanced property, the spatial discretization for both the cell average and the point values was specially designed, with the source term reformulated based on the target hydrostatic steady state solutions and a properly modified Lax-Friedrichs flux vector splitting for the flux term. The positivity-preserving property was achieved by applying a posteriori MOOD paradigm equipped with a structure-preserving first-order scheme. Theoretical analysis and numerical experiments have demonstrated our schemes' high-order accuracy, positivity preservation and robustness in capturing complex small features at or near the equilibrium flow. 

\setcounter{equation}{0}
\renewcommand\theequation{A.\arabic{equation}}
\renewcommand{\appendixname}{Appendix~\Alph{section}}
\begin{appendix}
\section{Useful properties }\label{a1}
Some important properties for the admissible state and the one-dimensional HLLC numerical flux are listed here, two-dimensional case is similar. All the lemmas can be found in \cite{wu2021uniformly} as well.
\begin{lemma}[\cite{wu2021uniformly}]\label{lemma:ad_s}
  For any constant $\lambda\geq 0$, $\delta\in \mathbb{R}$, $\boldsymbol{U}= (\rho,\rho\boldsymbol{u},E)^T\in \mathcal{G}$, and $\boldsymbol{a}\in \mathbb{R}^d$, if $\delta\frac{||a||}{\sqrt{2e}}\leq \lambda$, then $\lambda\boldsymbol{U}+\delta(0,\rho\boldsymbol{a},\rho\boldsymbol{u}\cdot\boldsymbol{a})^T \in \overline{\mathcal{G}}$.
\end{lemma}

\begin{lemma}[\cite{wu2021uniformly}]\label{lemma:ad_f}
  Let $\boldsymbol{H}(\boldsymbol{U})=[\boldsymbol{F}(\boldsymbol{U}),\boldsymbol{G}(\boldsymbol{U})]$ be the flux functions. For any $\boldsymbol{U}\in\mathcal{G}$ and unit vector $\boldsymbol{n}\in\mathbb{R}^d$, we have $\boldsymbol{U}-\lambda\boldsymbol{H}(\boldsymbol{U})\cdot\boldsymbol{n}\in\mathcal{G}$, for any $\lambda\in\mathbb{R}$ satisfying 
  $|\lambda|\varrho(\boldsymbol{U})\leq 1$.
\end{lemma}

\begin{lemma}[Contact property\cite{wu2021uniformly}]\label{lemma:contact}
  For any two states $\boldsymbol{U}_l = (\rho_l,0,p/(\gamma-1))^T$,  and $\boldsymbol{U}_r = (\rho_r,0,p/(\gamma-1))^T$, HLLC flux satisfies $\boldsymbol{F}^{hllc}(\boldsymbol{U}_l,\boldsymbol{U}_r) = (0,p,0)^T$.
\end{lemma}

\begin{lemma}[\cite{wu2021uniformly}]\label{lemma:cfl}
  For any parameters $\zeta_1,\zeta_2,\zeta_3,\zeta_4\geq 0$ and any admissible states $\boldsymbol{U}_l$, $\boldsymbol{U}_r \in \mathcal{G}$, 
  if $\lambda>0$ satisfies $\lambda\max_{\boldsymbol{U}\in\{\boldsymbol{U}_l, \boldsymbol{U}_r\}} \varrho(\boldsymbol{U})\leq \frac{1}{2}$, then we have 
  \begin{align*}
    &\zeta_2\boldsymbol{U}_r - \lambda(\boldsymbol{F}(\zeta_2\boldsymbol{U}_r) - \boldsymbol{F}^{hllc}(\zeta_1\boldsymbol{U}_l,\zeta_2\boldsymbol{U}_r)) \in \mathcal{G}, \\
    &\zeta_3\boldsymbol{U}_l - \lambda(\boldsymbol{F}^{hllc}(\zeta_3\boldsymbol{U}_l,\zeta_4\boldsymbol{U}_r)-\boldsymbol{F}(\zeta_3\boldsymbol{U}_l)) \in \mathcal{G}.
  \end{align*} 
\end{lemma}

\section{Proof of Theorem \ref{wb_1st_p} and \ref{pp_1st_p}}\label{a2}
\begin{proof}[Proof of the well-balanced property \ref{wb_1st_p}]
  For brevity, we prove only \eqref{1st_x} here, as the proofs for \eqref{1st_y} and \eqref{1st_c} proceed analogously. 
  
  At the equilibrium state, \eqref{steady_ini} is satisfied, and we have
  \begin{equation*}
    \frac{p^{e,*}_{i+1,j}}{ p^e_{i+\frac{1}{2},j}}{\boldsymbol{U}}_{i+\frac{1}{2},j} = \left(\frac{p^{e,*}_{i+1,j}}{ p^e_{i+\frac{1}{2},j}}\rho_{i+\frac{1}{2},j}, 0,0,\frac{p^{e,*}_{i+1,j}}{\gamma-1}\right)^T, \ \frac{p^{e,*}_{i+1,j}}{ p^e_{i+\frac{3}{2},j}}{\boldsymbol{U}}_{i+\frac{3}{2},j} = \left(\frac{p^{e,*}_{i+1,j}}{ p^e_{i+\frac{3}{2},j}}\rho_{i+\frac{3}{2},j}, 0,0,\frac{p^{e,*}_{i+1,j}}{\gamma-1}\right)^T.
  \end{equation*}
  According to the contract property Lemma \ref{lemma:contact}, the modified HLLC flux becomes $\widehat{\boldsymbol{F}}^L_{i+1,j} = (0,p^{e,*}_{i+1,j},0,0)^T$. Similarly, we have $\widehat{\boldsymbol{G}}^L_{i+\frac{1}{2},j+\frac{1}{2}} = (0,0,p^{e,*}_{i+\frac{1}{2},j+\frac{1}{2}},0)^T$. Owing to the fact $(\rho/\rho^e)_{i+\frac{1}{2},j} = 1$ in equilibrium, and the central difference approximation for the source term, we have
  \begin{equation*}
    \frac{1}{2}(p^{e}_{i+\frac{3}{2},j} - p^{e}_{i-\frac{1}{2},j}) = p^{e,*}_{i+1,j} - p^{e,*}_{i,j},\quad \frac{1}{2}(p^e_{i+\frac{1}{2},j+1} - p^e_{i+\frac{1}{2},j-1}) = p^{e,*}_{i+\frac{1}{2},j+\frac{1}{2}} - p^{e,*}_{i+\frac{1}{2},j-\frac{1}{2}}.
  \end{equation*}
  Thus the momentum equation in \eqref{1st_x} equals
  \begin{align*}
    \dfrac{d}{dt}\boldsymbol{U}_{i+\frac{1}{2},j}^{(2)} &= -\frac{1}{\Delta x}(p^{e,*}_{i+1,j} - p^{e,*}_{i,j}) + \frac{1}{2\Delta x}(p^{e}_{i+\frac{3}{2},j} - p^{e}_{i-\frac{1}{2},j}) \overset{\eqref{pe_1st}}{=}0, \\
    \dfrac{d}{dt}\boldsymbol{U}_{i+\frac{1}{2},j}^{(3)} &= -\frac{1}{\Delta y}(p^{e,*}_{i+\frac{1}{2},j+\frac{1}{2}} - p^{e,*}_{i+\frac{1}{2},j-\frac{1}{2}}) + \frac{1}{2\Delta y}(p^e_{i+\frac{1}{2},j+1} - p^e_{i+\frac{1}{2},j-1}) \overset{\eqref{pe_1st}}{=}0.
  \end{align*}
  Moreover, the well-balanced property for the mass and energy equations in \eqref{1st_x} follows directly, since the numerical flux and source term both reduce to zero. This completes the proof of Theorem \ref{wb_1st_p}.
\end{proof}

\begin{proof}[Proof of the positivity-preserving property \ref{pp_1st_p}]
  In fact, the proof follows essentially the same strategy as that in \cite{wu2021uniformly}.
  The scheme \eqref{1st_x} coupled with Euler forward time discretization yields
  \begin{small}
  \begin{align*}
    \boldsymbol{U}_{i+\frac{1}{2},j}^{n+1} = \boldsymbol{U}_{i+\frac{1}{2},j}&- 
    \frac{\Delta t}{\Delta x}(\widehat{\boldsymbol{F}}^L_{i+1,j} - \widehat{\boldsymbol{F}}^L_{i,j}) - \frac{\Delta t}{\Delta y}(\widehat{\boldsymbol{G}}^L_{i+\frac{1}{2},j+\frac{1}{2}} - \widehat{\boldsymbol{G}}^L_{i+\frac{1}{2},j-\frac{1}{2}})+\Delta t\boldsymbol{S}_{i+\frac{1}{2},j}\\
    =\boldsymbol{U}_{i+\frac{1}{2},j}&-
    \frac{\Delta t}{\Delta x}\frac{p^{e,*}_{i+1,j}-p^{e,*}_{i,j}}{p^e_{i+\frac{1}{2},j}}\boldsymbol{F}({\boldsymbol{U}}_{i+\frac{1}{2},j})-\frac{\Delta t}{\Delta y}\frac{p^{e,*}_{i+\frac{1}{2},j+\frac{1}{2}} - p^{e,*}_{i+\frac{1}{2},j-\frac{1}{2}}}{p^e_{i+\frac{1}{2},j}}\boldsymbol{G}({\boldsymbol{U}}_{i+\frac{1}{2},j})+\Delta t\boldsymbol{S}_{i+\frac{1}{2},j}\\
    &+\frac{\Delta t}{\Delta x}\left[\boldsymbol{F}\left(\frac{p^{e,*}_{i+1,j}}{ p^e_{i+\frac{1}{2},j}}{\boldsymbol{U}}_{i+\frac{1}{2},j}\right) - \boldsymbol{F}^{hllc}\left(\frac{p^{e,*}_{i+1,j}}{ p^e_{i+\frac{1}{2},j}}{\boldsymbol{U}}_{i+\frac{1}{2},j}, \frac{p^{e,*}_{i+1,j}}{ p^e_{i+\frac{3}{2},j}}{\boldsymbol{U}}_{i+\frac{3}{2},j} \right) \right]\\
    &- \frac{\Delta t}{\Delta x}\left[
    \boldsymbol{F}\left(\frac{p^{e,*}_{i,j}}{ p^e_{i+\frac{1}{2},j}}{\boldsymbol{U}}_{i+\frac{1}{2},j}\right) - 
    \boldsymbol{F}^{hllc}\left(\frac{p^{e,*}_{i,j}}{ p^e_{i-\frac{1}{2},j}}{\boldsymbol{U}}_{i-\frac{1}{2},j}, \frac{p^{e,*}_{i,j}}{ p^e_{i+\frac{1}{2},j}}{\boldsymbol{U}}_{i+\frac{1}{2},j} \right) \right]\\
    &+\frac{\Delta t}{\Delta y}\left[\boldsymbol{G}\left(\frac{p^{e,*}_{i+\frac{1}{2},j+\frac{1}{2}}}{p^e_{i+\frac{1}{2},j}}\boldsymbol{U}_{i+\frac{1}{2},j}\right) - \boldsymbol{G}^{hllc}\left(\frac{p^{e,*}_{i+\frac{1}{2},j+\frac{1}{2}}}{p^e_{i+\frac{1}{2},j}}\boldsymbol{U}_{i+\frac{1}{2},j}, \frac{p^{e,*}_{i+\frac{1}{2},j+\frac{1}{2}}}{p^e_{i+\frac{1}{2},j+1}}\boldsymbol{U}_{i+\frac{1}{2},j+1} \right)\right]\\
    &-\frac{\Delta t}{\Delta y}\left[ \boldsymbol{G}\left(\frac{p^{e,*}_{i+\frac{1}{2},j-\frac{1}{2}}}{p^e_{i+\frac{1}{2},j}}\boldsymbol{U}_{i+\frac{1}{2},j}\right) - \boldsymbol{G}^{hllc}\left(\frac{p^{e,*}_{i+\frac{1}{2},j-\frac{1}{2}}}{p^e_{i+\frac{1}{2},j-1}}\boldsymbol{U}_{i+\frac{1}{2},j-1}, \frac{p^{e,*}_{i+\frac{1}{2},j-\frac{1}{2}}}{p^e_{i+\frac{1}{2},j}}\boldsymbol{U}_{i+\frac{1}{2},j} \right)\right].
  \end{align*}
\end{small}
Here the superscript ${}^n$ is omitted from the right hand side.
Subsequently, we split $\boldsymbol{U}_{i+\frac{1}{2},j}^{n+1}$ into four parts, i.e. $\boldsymbol{U}_{i+\frac{1}{2},j}^{n+1} := \boldsymbol{W}_1+\boldsymbol{W}_2+\boldsymbol{W}_3+\boldsymbol{W}_4$, with
\begin{footnotesize}
\begin{align*}
  \boldsymbol{W}_1 &= \left[1-\Delta t\left(2\frac{\widehat\alpha_F}{p^{e}_{i+\frac{1}{2},j}}\left(\frac{p^{e,*}_{i+1,j}+p^{e,*}_{i,j}}{\Delta x}+\frac{p^{e,*}_{i+\frac{1}{2},j+\frac{1}{2}}+p^{e,*}_{i+\frac{1}{2},j-\frac{1}{2}}}{\Delta y}\right)+\widehat\alpha_S\right)\right]\boldsymbol{U}_{i+\frac{1}{2},j},\\
  \boldsymbol{W}_2 &=\Delta t \left[\frac{\widehat\alpha_F}{\Delta x}\frac{p^{e,*}_{i+1,j}}{ p^e_{i+\frac{1}{2},j}}\left({\boldsymbol{U}}_{i+\frac{1}{2},j} - \frac{1}{\widehat\alpha_F}\boldsymbol{F}({\boldsymbol{U}}_{i+\frac{1}{2},j})\right) + \frac{\widehat\alpha_F}{\Delta x}\frac{p^{e,*}_{i,j}}{ p^e_{i+\frac{1}{2},j}}\left({\boldsymbol{U}}_{i+\frac{1}{2},j} + \frac{1}{\widehat\alpha_F}\boldsymbol{F}({\boldsymbol{U}}_{i+\frac{1}{2},j})\right)\right.\\
  &\qquad+\left.\frac{\widehat\alpha_F}{\Delta y}\frac{p^{e,*}_{i+\frac{1}{2},j+\frac{1}{2}}}{p^e_{i+\frac{1}{2},j}}\left(\boldsymbol{U}_{i+\frac{1}{2},j}-\frac{1}{\widehat\alpha_F}\boldsymbol{G}({\boldsymbol{U}}_{i+\frac{1}{2},j})\right) + \frac{\widehat\alpha_F}{\Delta y}\frac{p^{e,*}_{i+\frac{1}{2},j-\frac{1}{2}}}{p^e_{i+\frac{1}{2},j}}\left(\boldsymbol{U}_{i+\frac{1}{2},j}+\frac{1}{\widehat\alpha_F}\boldsymbol{G}({\boldsymbol{U}}_{i+\frac{1}{2},j})\right)\right],\\
  \boldsymbol{W}_3 &=\Delta t \left(\widehat\alpha_S\boldsymbol{U}_{i+\frac{1}{2},j} + \boldsymbol{S}_{i+\frac{1}{2},j}\right),\\
  \boldsymbol{W}_4 &=\Delta t \left\{\frac{\widehat\alpha_F}{\Delta x}\left[\frac{p^{e,*}_{i+1,j}}{ p^e_{i+\frac{1}{2},j}}{\boldsymbol{U}}_{i+\frac{1}{2},j}-\frac{1}{\widehat\alpha_F}\left[\boldsymbol{F}^{hllc}\left(\frac{p^{e,*}_{i+1,j}}{ p^e_{i+\frac{1}{2},j}}{\boldsymbol{U}}_{i+\frac{1}{2},j}, \frac{p^{e,*}_{i+1,j}}{ p^e_{i+\frac{3}{2},j}}{\boldsymbol{U}}_{i+\frac{3}{2},j} \right)-\boldsymbol{F}\left(\frac{p^{e,*}_{i+1,j}}{ p^e_{i+\frac{1}{2},j}}{\boldsymbol{U}}_{i+\frac{1}{2},j}\right) \right]\right]\right.\\
  &\qquad+\frac{\widehat\alpha_F}{\Delta x}\left[\frac{p^{e,*}_{i,j}}{ p^e_{i+\frac{1}{2},j}}{\boldsymbol{U}}_{i+\frac{1}{2},j}-\frac{1}{\widehat\alpha_F}\left[\boldsymbol{F}\left(\frac{p^{e,*}_{i,j}}{ p^e_{i+\frac{1}{2},j}}{\boldsymbol{U}}_{i+\frac{1}{2},j}\right)-
  \boldsymbol{F}^{hllc}\left(\frac{p^{e,*}_{i,j}}{ p^e_{i-\frac{1}{2},j}}{\boldsymbol{U}}_{i-\frac{1}{2},j}, \frac{p^{e,*}_{i,j}}{ p^e_{i+\frac{1}{2},j}}{\boldsymbol{U}}_{i+\frac{1}{2},j} \right) \right]\right]\\
  &\qquad +\frac{\widehat\alpha_F}{\Delta y}\left[\frac{p^{e,*}_{i+\frac{1}{2},j+\frac{1}{2}}}{p^e_{i+\frac{1}{2},j}}\boldsymbol{U}_{i+\frac{1}{2},j} - \frac{1}{\widehat\alpha_F}\left[ \boldsymbol{G}^{hllc}\left(\frac{p^{e,*}_{i+\frac{1}{2},j+\frac{1}{2}}}{p^e_{i+\frac{1}{2},j}}\boldsymbol{U}_{i+\frac{1}{2},j}, \frac{p^{e,*}_{i+\frac{1}{2},j+\frac{1}{2}}}{p^e_{i+\frac{1}{2},j+1}}\boldsymbol{U}_{i+\frac{1}{2},j+1} \right)-\boldsymbol{G}\left(\frac{p^{e,*}_{i+\frac{1}{2},j+\frac{1}{2}}}{p^e_{i+\frac{1}{2},j}}\boldsymbol{U}_{i+\frac{1}{2},j}\right) \right]\right]\\
  &\qquad +\left.\frac{\widehat\alpha_F}{\Delta y}\left[\frac{p^{e,*}_{i+\frac{1}{2},j-\frac{1}{2}}}{p^e_{i+\frac{1}{2},j}}\boldsymbol{U}_{i+\frac{1}{2},j} - \frac{1}{\widehat\alpha_F}\left[ \boldsymbol{G}\left(\frac{p^{e,*}_{i+\frac{1}{2},j-\frac{1}{2}}}{p^e_{i+\frac{1}{2},j}}\boldsymbol{U}_{i+\frac{1}{2},j}\right)-\boldsymbol{G}^{hllc}\left(\frac{p^{e,*}_{i+\frac{1}{2},j-\frac{1}{2}}}{p^e_{i+\frac{1}{2},j-1}}\boldsymbol{U}_{i+\frac{1}{2},j-1}, \frac{p^{e,*}_{i+\frac{1}{2},j-\frac{1}{2}}}{p^e_{i+\frac{1}{2},j}}\boldsymbol{U}_{i+\frac{1}{2},j} \right)\right]\right]\right\}.
\end{align*}
\end{footnotesize}
Using the Lemma \ref{lemma:ad_s}, \ref{lemma:ad_f}, \ref{lemma:cfl} given in Appendix \ref{a1} and following the arguments in \cite{wu2021uniformly}, we can establish that $\boldsymbol{W}_i\in\mathcal{G}, i=1,\cdots,4$. Here we only prove the positivity-preserving property for \eqref{1st_x}, \eqref{1st_y} and \eqref{1st_c} can be treated in a similar manner.
\end{proof}

\end{appendix}

\section*{Declarations}

\subsection*{Data Availability}
The datasets generated during the current study are available from the corresponding author upon reasonable request. They support
our published claims and comply with field standards.

\subsection*{Competing of interest}
The authors declare that they have no known competing financial interests or personal relationships that could have appeared to
influence the work reported in this paper.

\subsection*{Funding}
The research of Christian Klingenberg was partially supported by the German Science Foundation (DFG) KL 566/24-1 as part of a joint project with the Spanish Agencia Estatal de Investigaci\'on (AEI).

\addcontentsline{toc}{section}{References}
\bibliographystyle{abbrv}
\normalem
\bibliography{shortreference}

\end{document}